\documentclass{article}

\usepackage{geometry}
\usepackage{graphicx} 
\usepackage{epstopdf} 
\usepackage[colorlinks, hypertexnames=false]{hyperref}
\hypersetup{linkcolor=blue, citecolor=red} 

\usepackage{bm}
\usepackage{amsmath} 
\usepackage{amssymb} 
\usepackage{stmaryrd} 
\usepackage{multirow}
\usepackage{booktabs}
\usepackage{subfigure}
\newcommand{\hong}[1]{{\color{black}{#1}}}

\usepackage{amsthm} 
\usepackage{mathtools} 

\usepackage[dvipsnames]{xcolor}
\numberwithin{equation}{section}
\usepackage{adjustbox}
\newtheorem{theorem}{Theorem}[section]
\newtheorem{lemma}{Lemma}[section]
\newtheorem{remark}{Remark}[section]

\graphicspath{{./Figs/}}
\renewcommand{\footnotesize}{\scriptsize}
\newcommand{\normmm}[1]{{\left\vert\kern-0.25ex\left\vert
		\kern-0.25ex\left\vert #1
		\right\vert\kern-0.25ex\right\vert\kern-0.25ex\right\vert}}

\usepackage{tikz}
\usetikzlibrary{math}
\usetikzlibrary{calc}

\newcommand{\jump}[1]{\ensuremath{\left\llbracket #1\right\rrbracket}}
\newcommand{\tdiv}{\mathrm{div}}
\newcommand{\tcurl}{\mathrm{curl}}

\newcommand{\Red}[1]{\textcolor{black}{#1}}
\newcommand{\Blue}[1]{\textcolor{black}{#1}}
\newcommand{\vertiii}[1]{{\left\vert\kern-0.25ex\left\vert\kern-0.25ex\left\vert #1 
		\right\vert\kern-0.25ex\right\vert\kern-0.25ex\right\vert}}

\begin{document}

	\title{A novel multipoint stress control volume method for linear elasticity on quadrilateral grids}
\author{Shubin Fu\footnotemark[1] \quad and \quad Lina Zhao\footnotemark[2]}
\renewcommand{\thefootnote}{\fnsymbol{footnote}}
\footnotetext[1]{Eastern Institute for Advanced Study, Eastern Institute of Technology, Zhejiang 315200,  PR China. }
\footnotetext[2]{Lina Zhao is the corresponding author. Department of Mathematics, City University
	of Hong Kong, Kowloon Tong, Hong Kong SAR, China. ({linazha@cityu.edu.hk}). 
The research of Lina Zhao is partially supported by the Research Grants Council of the Hong Kong Special Administrative Region, China. (Project No. CityU 21309522).}

\maketitle
\begin{abstract}
	
	In this paper, we develop a novel control volume method that is locally conservative and locking-free for linear elasticity problem on quadrilateral grids. The symmetry of stress is weakly imposed through the introduction of a Lagrange multiplier. As such, the method involves three unknowns: stress,  displacement and rotation. To ensure the well-posedness of the scheme, a pair of carefully defined finite element spaces is used for the stress, displacement and rotation such that the inf-sup condition holds. An appealing feature of the method is that piecewise constant functions are used for the approximations of stress, displacement and rotation, which greatly simplifies the implementation. In particular, the stress space is defined delicately such that the stress bilinear form is localized around each vertex, which allows for the local elimination of the stress, resulting in a cell-centered system. By choosing different definitions of the space for rotation, we develop two variants of the method. In particular, the first method uses a constant function for rotation over the interaction region, which allows for  further elimination and results in a cell-centered system involving displacement only. A rigorous error analysis is performed for the proposed scheme. We show the optimal convergence for $L^2$-error of the stress and rotation. Moreover, we can also prove the superconvergence for $L^2$-error of displacement. Extensive numerical simulations indicate that our method is efficient and accurate, and can handle problems with discontinuous coefficients. 
	
\end{abstract}



%

\section{Introduction}
Mixed finite element methods in stress-displacement formulation for linear elasticity is popular in solid mechanics since they avoid locking and provide a direct approximation to the stress that is the primary physical interest. Numerous methods have been developed in the context of strong stress symmetry \cite{Arnold84,Arnold02,Gopalakrishnan11,PechsteinTDNNS18} and weak stress symmetry \cite{Arnold84peer,arnold2007mixed,Boffi09,Zhao20}. From the computational point of view, the introduction of the stress element will lead to a saddle-point system and result in additional computational costs, therefore, many of these methods could suffer from extra computational costs. The common approaches to overcome this issue lie in hybridization and reduction to cell-centered system. The former has been successfully applied in the context of nonconforming mixed finite element (MFE) methods \cite{Arnold14} and hybridizable discontinuous Galerkin methods \cite{Qiu18}. The latter is to design suitable strategies to treat the stress part such that the local elimination can be applied to obtain a reduced system. In this context, we mention in particular the multipoint stress mixed finite element method \cite{Ambartsumyan20simplicial,Ambartsumyan21}. Therein, MFE spaces with the lowest order Brezzi-Douglas-Marini ($\text{BDM}_1$) degrees of freedom is used for the stress and piecewise constant approximation is used for the displacement. In addition, the vertex quadrature rule is used for the computation of the stress bilinear form, which localizes the stress degrees of freedom. As such, the mass matrix for stress is block-diagonal, and thus allows local elimination of the stress. The resulting system is symmetric and positive definite, which enhances the computational efficiency. This method is motivated by the multipoint flux mixed finite element (MFMFE) method \cite{Wheeler06,Ingram10,wheeler2012multipoint} for Darcy flow that is closely related to the multipoint flux approximation (MPFA) method \cite{Aavatsmark98,Edwards98,Aavatsmark02,Edwards02,Agelas08}. The MFMFE method invokes $\text{BDM}_1$ on simplicial and quadrilateral grids. As an alternative, a MFEM based on broken Raviart-Thomas velocity space is proposed in \cite{klausen2006robust,Klausen06}. All these methods share the similar idea that a vertex quadrature rule is applied for the computation of the mass matrix, which results in a block-diagonal mass matrix. Therefore, the flux can be locally eliminated, which leads to a cell-centered pressure system, rendering the method computationally attractive.

In this paper we aim to develop a new method that only uses piecewise constant approximations for the involved unknowns, that can be further reduced to a symmetric positive definite system. Our method is closely related to the multipoint stress mixed finite element method proposed in \cite{Ambartsumyan20simplicial,Ambartsumyan21}, but with much simpler  construction and implementation in the sense that no special quadrature rule is needed. The devising of the stress space is more subtle as it needs to be carefully balanced with the displacement such that the inf-sup condition holds. To this end, we divide the quadrilateral element into four smaller quadrilaterals by connecting the interior points to the midpoint of each edge and then define the stress space as a constant function over each smaller quadrilateral, and at the same time it is normal continuous over the edges lying on the quadrilaterals, but no continuity is imposed for the new edges generated by the subdivision. As a consequence, the bilinear form associated with the stress is localized around each vertex, which resembles the vertex quadrature rule used in multipoint stress mixed finite element method. The major difference is that we use piecewise constant function to avoid the quadrature rule. Then we develop two variants of the method by choosing different spaces for the rotation. In the first method, we let the rotation be a constant function over each interaction region formed by the four smaller quadrilaterals sharing the original vertex. Then the system can be further reduced to achieve a symmetric and positive definite cell-centered system for the displacement only. In the second method, we choose the rotation to be a constant function over each quadrilateral. In this case, the rotation  can not be further eliminated, but this variant can handle discontinuous coefficients accurately.  In particular, a simple element-wise reconstruction will lead to $H(\tdiv;\Omega)$-conforming stress. Our method can also be viewed as a mixed-type finite volume method, but we can carry out the error analysis in the framework of finite element methods. This can simplify the analysis as it is well-known it is not always easy to analyze finite volume method. \Red{We emphasize that the current approach shares a conceptual similarity with the method proposed in \cite{Zhao20}, particularly in the use of primal and dual grids to ensure method stability. The first method introduced in this paper, after the elimination of stress and rotation variables, results in a smaller system compared to that of \cite{Zhao19}, potentially leading to faster computations. The second method presented here incorporates displacement and rotation in the final system, with its size determined by the number of primal elements. In contrast, the size of the final system in the method proposed by \cite{Zhao19} is dependent on the number of primal edges. Moreover, the method from \cite{Zhao19} is capable of handling highly distorted grids, however, the current methods require the $\mathcal{O}(h^2)$ mesh condition.}

We perform a rigorous error analysis for the proposed scheme. The unique solvability of the solution is established. Our method also inherits the locking-free property of standard mixed finite element method. Besides, we also prove the superconvergence of $L^2$-error of displacement under the assumption that the solution is smooth enough. 
   In summary, our method owns several appealing features, which include: (1) the method enjoys local conservation, and a simple reconstruction will yield $H(\tdiv;\Omega)$-conforming stress; (2) it is locking-free, which allows to handle problems with nearly incompressible materials; (3) piecewise constant functions are used for all the unknowns, which makes the implementation easy; (4) the stress bilinear form is block-diagonal, which allows local elimination of stress and reduction to a cell-centered system. All these salient features make our method a good candidate for handling problems arising from practical applications. 
   
   Extensive numerical simulations including three-dimensional tests are carried out. 
 We observe that both methods deliver at least first-order convergence for all variables in either structured or smooth unstructured grids, 
 superconvergence of displacements and rotation are observed. Moreover, by taking average of stress tensor inside each element also 
 leads to superconvergence in  structured and some types of smooth unstructured grids. Besides,  our methods are observed to be much more accurate than $\text{BDM}_1$ with mass lumping strategies, locking-free properties are found as expected and the ability to solve challenging 
 highly heterogeneous media are demonstrated.

The major contributions of the paper can be summarized as follows:
\begin{enumerate}
	\item We design a locally conservative and locking-free method in stress-displacement formulation with piecewise constant approximations that allows reduction to a cell-centered system. 
	\item We prove the convergence error estimates for all the involved variables. Moreover, the superconvergence for displacement can be proved under the assumption that the solution is smooth smooth.
	\item We perform extensive numerical simulations including three-dimensional tests to demonstrate the performance of the scheme. Furthermore, the numerical results indicate that our method has superior performance in terms of accuracy and efficiency. 
\end{enumerate}
The rest of the paper is organized as follows. In the next section, 
we introduce the model problem and present the corresponding weak formulation. In section~\ref{sec:method}, we describe the proposed methods. The comprehensive error analysis and the superconvergence for $L^2$-error of displacement are performed in section~\ref{sec:analysis}. \Blue{The extension of the proposed scheme to Darcy equations and the Stokes equations are given in section~\ref{sec:extension}.} In section~\ref{sec:numerical}, extensive numerical simulations including three-dimensional tests are carried out to demonstrate the capability, efficiency and superior performance of the proposed methods.

\section{The mathematical model}
\label{sec:main}

Let $\Omega\subset \mathbb{R}^d,d=2,3$ be a polygonal domain. We consider the following model problem 
\begin{alignat}{2}
	\mathcal{A}\underline{\sigma}&= \varepsilon(\bm{u}) &&\quad \mbox{in}\;\Omega,\label{eq:elasticity1}\\
	-\text{div}\,\underline{\sigma}&=\bm{f} &&\quad \mbox{in}\;\Omega,\label{eq:elasticity2}
\end{alignat}
where $\bm{u}=\bm{0}$ on $\partial \Omega$ and $\mathcal{A}$ is the compliance
tensor determined by material parameters of the elastic medium. In a homogeneous isotropic elastic medium, $\mathcal{A}$ has the form
\begin{align*}
	\mathcal{A}\underline{G}=\frac{1}{2\mu}(\underline{G}-\frac{\lambda}{d\lambda+2\mu}\mbox{tr}(\underline{G})\underline{I}).
\end{align*}
Here $\lambda$ and $\mu$ are positive, called the Lam\'{e} parameters, $\mbox{tr}(\underline{G})$ is the trace of function $\underline{G}$, and $\underline{I}$ is the identity matrix. We define $\|\underline{G}\|_{\mathcal{A}}^2=(\mathcal{A}\underline{G},\underline{G})$. 
Now we introduce some notation that will be used throughout the paper.  The norms and seminorms of the Sobolev spaces $W^{k,p}(\mathcal{O}), k\in \mathbb{R},p>0$ are denoted by $\|\cdot\|_{W^{k,p}(\mathcal{O})}$ and $|\cdot|_{W^{k,p}(\mathcal{O})}$, respectively, where $\mathcal{O}\subset \Omega,\Omega\subset \mathbb{R}^d$. When $\Blue{p=2}$, we use
$H^r(\mathcal{O})$ to represent the corresponding space. The spaces of vector- and matrix-valued functions
with all the components in $H^r(\mathcal{O})$ will be respectively denoted as $\bm{H}^r(\mathcal{O})$ and $\underline{H}^r(\mathcal{O})$. 
 The corresponding norm and semi-norm are respectively denoted as $\|\cdot\|_{H^r(\mathcal{O})}$ and $|\cdot|_{H^r(\mathcal{O})}$. 
We use $(\cdot,\cdot)_D$ to represent the standard $L^2$-inner product over $D\subset \mathbb{R}^d$ and when $D$ coincides with $\Omega$, the subscript will be omitted. The corresponding norm is denoted as $\|\cdot\|_{L^2(D)}$. 
We define $\bm{H}(\tdiv;\Omega):=\{\bm{v}\in \bm{L}^2(\Omega), \nabla\cdot \bm{v}\in L^2(\Omega)\}$, which is equipped with the norm $\|\bm{v}\|_{\tdiv}^2 = \|\bm{v}\|_{L^2(\Omega)}^2+\|\nabla\cdot \bm{v}\|_{L^2(\Omega)}^2$.
In addition, we use  $\underline{H}(\tdiv;\Omega)$ to represent the tensor field where each row belongs to $\bm{H}(\tdiv;\Omega)$. 

We define the following spaces:
\begin{align*}
	\underline{\Sigma}:=\{\underline{\tau}\in \underline{H}(\tdiv;\Omega)\},\quad
	\bm{U}:=\bm{L}^2(\Omega),\quad \Gamma:=L^2(\Omega).
\end{align*}
Moreover, we define the constant tensor
\begin{align*}\underline{\delta}=
	\begin{pmatrix}
		0&-1\\
		1 & 0
	\end{pmatrix}.
	\end{align*}
Introducing the Lagrange multiplier $\gamma =\text{rot} ( \bm{u})$, where $\text{rot}(\bm{u})=\frac{1}{2}(-\frac{\partial u_1}{\partial y}+\frac{\partial u_2}{\partial x})$ for $\bm{u}=(u_1,u_2)^T$, we have $\mathcal{A}\underline{\sigma} = \nabla \bm{u}-\gamma \underline{\delta}$. Then we can propose the 
following weak formulation when $d=2$: Find $(\underline{\sigma},\bm{u},\gamma)\in \underline{\Sigma}\times \bm{U}\times \Gamma$ such that
\begin{alignat*}{2}
	(\mathcal{A} \underline{\sigma}, \underline{w})+(\bm{u},\tdiv \underline{w})+(\text{as}(\underline{w}), \gamma)&=0&&\quad \forall \underline{w}\in \underline{\Sigma},\\
-	(\tdiv \underline{\sigma},\bm{v})&=(\bm{f},\bm{v})&&\quad \forall \bm{v}\in \bm{U},\\
	(\text{as}(\underline{\sigma}),\xi)&=0&&\quad \forall \xi\in \Gamma,
\end{alignat*}
where
$
	\text{as}(\underline{w})=\underline{w}_{12}-\underline{w}_{21}.
$
When $d=3$, we should define 
\begin{align*}
\underline{\Gamma}:=	\{\underline{w}\in \underline{L}^2(\Omega); \underline{w}=-\underline{w}^T\},\quad 
	\text{as}(\underline{w})=(\underline{w}_{32}-\underline{w}_{23},\underline{w}_{31}-\underline{w}_{13},\underline{w}_{21}-\underline{w}_{12} )^T.
\end{align*}
Moreover,
\begin{align*} \mathbb{E}(p)=
	\begin{pmatrix}
		0& -p_3&p_2\\p_3& 0 &-p_1\\
		-p_2& p_1& 0 
	\end{pmatrix}\quad \mbox{for}\;p\in \mathbb{R}^3.
\end{align*}
Then the weak formulation can be written as follows: 
Find $(\underline{\sigma},\bm{u},\underline{\gamma})\in \underline{\Sigma}\times \bm{U}\times\underline{ \Gamma}$ such that
\begin{alignat*}{2}
	(\mathcal{A} \underline{\sigma}, \underline{w})+(\bm{u},\tdiv \underline{w})+(\underline{w}, \gamma)&=0&&\quad \forall \underline{w}\in \underline{\Sigma},\\
	(\tdiv \underline{\sigma},\bm{v})&=(\bm{f},\bm{v})&&\quad \forall \bm{v}\in \bm{U},\\
	(\underline{\sigma},\underline{\xi})&=0&&\quad \forall \underline{\xi}\in \underline{\Gamma}.
\end{alignat*}
Here for any $\underline{w}\in \underline{L}^2(\Omega)$ and $\underline{\xi}\in \underline{\Gamma}$, we have
\begin{align*}
	(\underline{w},\underline{\xi}) = (\text{as}(\underline{w}), \mathbb{E}^{-1}(\underline{\xi})).
\end{align*}
The well-posedness of the weak formulation can be found in \cite{Arnold84peer,arnold2007mixed}, which is omitted here for simplicity.

\section{Description of the new scheme}\label{sec:method}
\begin{figure}[t]
	\begin{center}
		\begin{tikzpicture}[scale=1.6]
			\coordinate (O) at (0,0);
			\coordinate (A) at (1,0);
			\coordinate (B) at (1,1);
			\coordinate (C) at (0,1);
			\coordinate (D) at (-1,1);
			\coordinate (E) at (-1.2,0);
			\coordinate (F) at (-0.8,-1);
			\coordinate (G) at (-0.2,-1);
			\coordinate (H) at (1,-1);
			
			\coordinate (AC1) at ($0.25*(O)+0.25*(A)+0.25*(B)+0.25*(C)$);
			\coordinate (AC2) at ($0.25*(C)+0.25*(D)+0.25*(E)+0.25*(O)$);
			\coordinate (AC3) at ($0.25*(E)+0.25*(F)+0.25*(G)+0.25*(O)$);
			\coordinate (AC4) at ($0.25*(G)+0.25*(H)+0.25*(A)+0.25*(O)$);

			\coordinate (LC1) at ($0.5*(A)+0.5*(B)$);
			\coordinate (LC2) at ($0.5*(B)+0.5*(C)$);
			\coordinate (LC3) at ($0.5*(C)+0.5*(D)$);
			\coordinate (LC4) at ($0.5*(D)+0.5*(E)$);
			\coordinate (LC5) at ($0.5*(E)+0.5*(F)$);
			\coordinate (LC6) at ($0.5*(F)+0.5*(G)$);
			\coordinate (LC7) at ($0.5*(G)+0.5*(H)$);
			\coordinate (LC8) at ($0.5*(H)+0.5*(A)$);
			
			\coordinate (LLC1) at ($0.5*(O)+0.5*(A)$);
			\coordinate (LLC2) at ($0.5*(O)+0.5*(C)$);
			\coordinate (LLC3) at ($0.5*(O)+0.5*(E)$);
			\coordinate (LLC4) at ($0.5*(O)+0.5*(G)$);

			\draw[ultra thick] (A) -- (B) -- (C) -- (D) -- (E) -- (F) -- (G) -- (H) -- cycle;
			\draw[ultra thick] (A) -- (O) -- (E);
			\draw[ultra thick] (C) -- (O) -- (G);
			
			\draw[dashed, red, thick] (LC1) -- (AC1) -- (LLC2) -- (AC2) -- (LC4);
			\draw[dashed, red, thick] (LC3) -- (AC2) -- (LLC3) -- (AC3) -- (LC6);
			\draw[dashed, red, thick] (LC5) -- (AC3) -- (LLC4) -- (AC4) -- (LC8);
			\draw[dashed, red, thick] (LC7) -- (AC4) -- (LLC1) -- (AC1) -- (LC2);
		\end{tikzpicture}
	\hskip 1cm
	\begin{tikzpicture}[scale=1.4]
		\coordinate (O) at (0,0);
		\coordinate (A) at (1,0);
		\coordinate (B) at (1,1);
		\coordinate (C) at (0,1);
		\coordinate (D) at (-1,1);
		\coordinate (E) at (-1,0);
		\coordinate (F) at (-1,-1);
		\coordinate (G) at (0,-1);
		\coordinate (H) at (1,-1);
		
		\draw[fill, black] (B) circle (1pt);
		\draw[fill, black] (D) circle (1pt);
		\draw[fill, black] (F) circle (1pt);
		\draw[fill, black] (H) circle (1pt);
		
		\draw[ultra thick, dashed] (A) -- (B) -- (C) -- (D) -- (E) -- (F) -- (G) -- (H) -- cycle;
		\draw[ultra thick] (A) -- (O) -- (E);
		\draw[ultra thick] (C) -- (O) -- (G);
		
		\coordinate (LLC1) at ($0.5*(O)+0.5*(A)$);
		\coordinate (LLC2) at ($0.5*(O)+0.5*(C)$);
		\coordinate (LLC3) at ($0.5*(O)+0.5*(E)$);
		\coordinate (LLC4) at ($0.5*(O)+0.5*(G)$);
		
		\draw[violet,->,thick] ($(LLC1)+(-0.2,-0.25)$) -- ($(LLC1)+(-0.2,0.25)$);
		\draw[violet,->,thick] ($(LLC2)+(0.25,-0.2)$) -- ($(LLC2)+(-0.25,-0.2)$);
		\draw[violet,->,thick] ($(LLC3)+(0.2,0.25)$) -- ($(LLC3)+(0.2,-0.25)$);
		\draw[violet,->,thick] ($(LLC4)+(-0.25,0.2)$) -- ($(LLC4)+(0.25,0.2)$);
		
		\draw[fill, cyan] (LLC1) circle (1pt);
		\draw[fill, cyan] (LLC2) circle (1pt);
		\draw[fill, cyan] (LLC3) circle (1pt);
		\draw[fill, cyan] (LLC4) circle (1pt);

		\node[red] at (AC1) {$E_3$};
		\node[red] at (AC2) {$E_4$};
		\node[red] at (AC3) {$E_1$};
		\node[red] at (AC4) {$E_2$};
		
		\node[above right] at (B) {$u_3$};
		\node[above left] at (D) {$u_4$};
		\node[below left] at (F) {$u_1$};
		\node[below right] at (H) {$u_2$};
		
		\node[above] at (LLC1) {$\underline{\sigma}_2$};
		\node[left] at (LLC2) {$\underline{\sigma}_3$};
		\node[below] at (LLC3){$\underline{\sigma}_4$};
		\node[right] at (LLC4) {$\underline{\sigma}_1$};
		
		\node[below left] at (A) {$e_2$};
		\node[below right] at (C) {$e_3$};
		\node[above right] at (E) {$e_4$};
		\node[above left] at (G) {$e_1$};
	\end{tikzpicture}
\end{center}
	\caption{The macro-element cell is divided into four subcells (left); the black solid lines represent the primal edges and the red dashed lines represent the dual edges (inner subcell edges). The interaction region is formed by the union of four subcells $E_i, i=1,\cdots,4$ sharing the common vertex (right). The filled dots denote the cell displacement $\{\bm{u}_i\}$ and the blue dots denote the stress $\{\underline{\sigma}_i\}$.}
	\label{fig:dof-MPFA}
\end{figure}
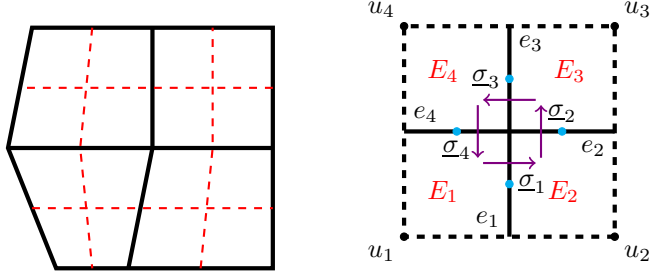

In this section, we will introduce the discrete formulation and state the main results of this paper. To simplify the presentation, we limit to $d=2$ for the construction and analysis of the proposed scheme. The proposed scheme can be extended to $d=3$ naturally. Moreover, we also carry out numerical simulations for $d=3$. To begin, 
we  introduce some notation that will be used throughout the paper. We let $\mathcal{T}_M$ represent the partition of the domain $\Omega$ into quadrilateral meshes. Each element $M\in \mathcal{T}_M$ is considered as the macro-element. The union of all the edges generated in this partition is denoted by $\mathcal{F}_{pr}$. For each macro-element $M$, we choose one interior point and connect it to the midpoints of the edges of $M$, then $M$ is decomposed into the union of four subcells.  Owing to this subdivision, each edge $e\in \mathcal{F}_{pr}$ is divided into two equal half edges, and the union of all the half edges is denoted as $\mathcal{F}_{pr}^{\frac{1}{2}}$; see Figure~\ref{fig:dof-MPFA} for an illustration. Furthermore, four inner subcell edges are generated in each macro-element, and the union of all the new edges lying inside of the macro-element is denoted as $\mathcal{F}_{dl}$. We let $\mathcal{F}_h:=\mathcal{F}_{pr}^{\frac{1}{2}}\cup \mathcal{F}_{dl}$. The union of subcells  is denoted as $\mathcal{T}_h$. In addition, we use $D$ to represent the interaction region, which is formed by the four subcells sharing the common vertex; see Figure~\ref{fig:dof-MPFA}. The union of all the interaction regions is denoted by $\mathcal{T}_D$.  We let $h_D$ represent the diameter of the element $D$, where $D$ could be the macro-element or the subcell and we let $h:=\max \{h_D\}$. Moreover, we use $h_e$ to represent the length of edge $e$. For on boundary facet $e$, we use $n_e$ to represent the unit normal vector of  pointing outside of $\Omega$. For an interior edge, we fix $\bm{n}_e$ as one of the two possible unit vectors. When there is no ambiguity, we use $\bm{n}$ to simplify the notation. Let $k\geq 0$ represent the polynomial order. We use $P_k$ to represent the  polynomial function whose order is less than or equal to $k$ and use $Q_k$ to represent the space of polynomials of order at most $k$ in each variable. 
Let $q$, $\bm{v}$ and $\underline{\omega}$ be scalar-, vector- and matrix-valued functions, respectively. For any two adjacent elements $E^+$ and $E^-$ sharing the common facet $e$, i.e., $e=\partial E^+\cap \partial E^-$, the jumps of $q$, $\bm{v}$ and $\underline{\omega}$ are given by
\begin{align*}
	\jump{q}:=q_{|E^+}-q_{|E^-},\quad \jump{\bm{v}}:=\bm{v}_{|E^+}-\bm{v}_{|E^-},\quad\jump{\underline{\omega}}:=\underline{\omega}_{|E^+}-\underline{\omega}_{|E^-}.
\end{align*}
On a boundary facet, we set  $\jump{q}:=q,\jump{\bm{v}}:=\bm{v},\jump{\underline{\omega}}:=\underline{\omega}$.

%
%
%
%
\subsection{Preliminaries}

\begin{itemize}
	\item(\textbf{Method 1}). We define the following spaces, which are important to define our multipoint stress control volume method:
	\begin{align*}
		\bm{\Sigma}_h^*:&=\{\bm{w}_h\in \bm{P}_0(E),\;\forall E\in \mathcal{T}_h, \jump{\bm{w}_h\cdot\bm{n}}_{|e}=0,\;\forall e\in \mathcal{F}_{pr}^{\frac{1}{2}}\backslash \partial \Omega\},\\
		\bm{U}_h:&=[U_h]^2=\{v_h\in P_0(M),\forall M\in \mathcal{T}_M\},\\
		\Gamma_h^1:&=\{\mu\in P_0(D),\forall D\in \mathcal{T}_D\}.
	\end{align*}
The right panel of	Figure~\ref{fig:dof-MPFA} displays the degrees of freedom of $\bm{\Sigma}_h^*$ and $U_h$ over each interaction region, which also indicates the stress is localized around each interaction region and can be locally eliminated, which inherits the similar idea to that of multipoint stress mixed finite element method \cite{Wheeler06}. 
	\item(\textbf{Method 2}). The stress and displacement spaces of Method 2 are the same to that of Method 1, the only difference lies in the definition of the rotation space. To this end, we define the following finite dimensional space:
	\begin{align*}
		\Gamma_h^2:&=\{\mu_h\in P_0(M),\forall M\in \mathcal{T}_M\}.
	\end{align*}
	We remark that the differences between these two method are the choice of the space for the rotation. In Method 1, we can apply local elimination for both stress and rotation. In Method 2, we can only apply local elimination for stress. More details can be found in section~\ref{sec:cell-center}. 
\end{itemize}

The basis functions for $\bm{\Sigma}_h^*$ over each interaction region (cf. Figure~\ref{fig:dof-MPFA}) are defined as
\begin{equation*}\bm{v}_1=
	\begin{cases}
		\frac{1}{\bm{t}_4\times \bm{t}_1}(t_4^1, t_4^2), \quad (x,y)\in E_1\\
		\frac{1}{\bm{t}_2\times \bm{t}_1}(t_2^1, t_2^2), \quad (x,y)\in E_2
	\end{cases},\quad\bm{v}_2=
	\begin{cases}
		\frac{1}{\bm{t}_1\times \bm{t}_2}(t_1^1, t_1^2), \quad (x,y)\in E_2\\
		\frac{1}{\bm{t}_3\times \bm{t}_2}(t_3^1, t_3^2), \quad (x,y)\in E_3
	\end{cases},
\end{equation*}

\begin{equation*}\bm{v}_3=
	\begin{cases}
		\frac{1}{\bm{t}_2\times \bm{t}_3}(t_2^1, t_2^2), \quad (x,y)\in E_3\\
		\frac{1}{\bm{t}_4\times \bm{t}_3}(t_4^1, t_4^2), \quad (x,y)\in E_4
	\end{cases},
\quad \bm{v}_4=
	\begin{cases}
		\frac{1}{\bm{t}_3\times \bm{t}_4}(t_3^1, t_3^2), \quad (x,y)\in E_4\\
		\frac{1}{\bm{t}_1\times \bm{t}_4}(t_1^1, t_1^2), \quad (x,y)\in E_1
	\end{cases},
\end{equation*}
where $\bm{t}_i:=(t_i^1,t_i^2)^T,i=1,\cdots,4$ represents the unit tangential vector of $e_i$. \Red{Here $\bm{a}\times \bm{b}:= a_1b_2-a_2b_1$ for $\bm{a}=(a_1,a_2)^T$ and $\bm{b}=(b_1,b_2)^T$.}

When the quadrilateral grid is reduced to the rectangular gird, then we can get the following simplified basis functions
\begin{align*}\bm{v}_1&=
	(1, 0)\quad   E_1\cup E_2,\quad \bm{v}_2=(0,1)\quad  E_2\cup E_3,\\
	\bm{v}_3&=(-1,0)\quad E_3\cup E_4\quad \bm{v}_4=(0,-1) \quad E_4\cup E_1.
\end{align*}
One can observe from the above definitions that the basis functions associated with $\underline{\Sigma}_h$ are rather easy to construct. Since we only use piecewise constant functions, the implementation is pretty simple, moreover, direct implementation on the physical element is enough. Compared to the methods proposed in \cite{Ambartsumyan21}, our method is easier to construct, and no special quadrature rule is needed. In particular, our numerical simulations indicate that the accuracy of our method is quite good. 

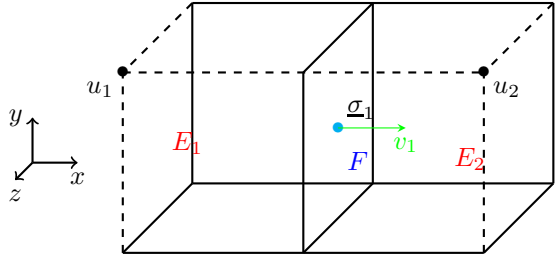
\begin{figure}
	\def\OverLen{0}
	\def\OverLenb{1.0}
	\def\LENGTH{4.0}
	\def\Opacity{0.25}
	\centering
	\begin{tikzpicture}[scale=.6]
		\coordinate (CoodOrigin) at ({-\LENGTH / 2}, {\LENGTH / 2}, \LENGTH);
		\draw[->,black, thick] (CoodOrigin) -- ++(\OverLenb, 0, 0) node[below] {$x$};
		\draw[->,black,thick] (CoodOrigin) -- ++(0, \OverLenb, 0) node[left] {$y$};
		\draw[->,black,thick] (CoodOrigin) -- ++(0, 0, \OverLenb) node[below] {$z$};
		
		\draw[thick] (0, 0, 0) -- ++(8, 0, 0);
		\draw[thick] (0, 4, 0) -- ++(8, 0, 0);
		\draw[thick] (0, {-0.5*\OverLen}, 0) -- ++(0, {{\LENGTH+\OverLen}}, 0);
		\draw[thick] (8, 0, 0) -- ++(0, 4, 0);
		\draw[thick] (0, 0,4) -- ++(8, 0, 0);
		\draw[thick, dashed] (0, 4, 4) -- ++(8, 0, 0);
		\draw[thick, dashed] (0, {-0.5*\OverLen}, \LENGTH) -- ++(0, {{\LENGTH+\OverLen}}, 0);
		
		\draw[thick] (\LENGTH, {-0.5*\OverLen}, \LENGTH) -- ++(0, {{\LENGTH+\OverLen}}, 0);
		\draw[thick] (\LENGTH, 0, {-0.5*\OverLen}) -- ++(0, 0, {\LENGTH+\OverLen});
		\draw[thick] (\LENGTH, \LENGTH, {-0.5*\OverLen}) -- ++(0, 0, {\LENGTH+\OverLen});
		\draw[thick, dashed] (0, \LENGTH, {-0.5*\OverLen}) -- ++(0, 0, {\LENGTH+\OverLen});
		\draw[thick] (0, 0, {-0.5*\OverLen}) -- ++(0, 0, {\LENGTH+\OverLen});
		
		\draw[thick] (8, 0, {-0.5*\OverLen}) -- ++(0, 0, {\LENGTH+\OverLen});
		\draw[thick, dashed] (8, 4 {-0.5*\OverLen}) -- ++(0, 0, {\LENGTH+\OverLen});
		\draw[thick, dashed] (8, {-0.5*\OverLen}, \LENGTH) -- ++(0, {{\LENGTH+\OverLen}}, 0);
		\draw[thick] (4, {-0.5*\OverLen}, 0) -- ++(0, {{\LENGTH+\OverLen}}, 0);
		
		\node[red, below right] at ({0.2*\LENGTH}, {0.7*\LENGTH}, {0.95*\LENGTH}) {$E_1$};
		\node[red, below right] at ({1.65*\LENGTH}, {0.5*\LENGTH}, {0.65*\LENGTH}) {$E_2$};
		\node[blue, below right] at ({\LENGTH}, {0.42*\LENGTH}, {0.5*\LENGTH}) {$F$};
		
		\node[above right] at  ({\LENGTH}, {0.5*\LENGTH}, {0.5*\LENGTH}) {$\underline{\sigma}_1$};
		\node[cyan] at  ({\LENGTH}, {0.5*\LENGTH}, {0.5*\LENGTH}) {$\bullet$};
		
		\node[below left] at (0,4,4) {$u_1$};
		\node[below right] at (8,4,4) {$u_2$};
		\node[black] at (0,4,4) {$\bullet$};
		\node[black] at (8,4,4) {$\bullet$};

		\draw[green,-stealth] (\LENGTH, {0.5*\LENGTH}, {0.5*\LENGTH}) -- ++(1.5, 0, 0) node[below]  {${v}_{1}$};
		
	\end{tikzpicture}
	
	\caption{Two subcells $E_1$ and $E_2$ associated with a common face $F$ in 3D, $u_1$ and $u_2$ are  the centers of two macro elements, each 
		black dot denotes a cell displacement, cyan dot denotes the stress $\underline{\sigma}_1$, the dashed edges are edges in dual grid, the solid edges are edges in macro-element.}
	\label{fig:3d}
\end{figure}

In 3D, we only consider the partition of grid as rectangular cuboid and the basis functions are defined in a similar way to the above. For example, let  $E_1$ 
and $E_2$ are two subcells in 3D associated with a face $F$ (see Figure \ref{fig:3d}), then we can define the basis function as
\begin{equation*}
	\bm{v}_1=(1,0,0),\quad E_1\cup E_2.
\end{equation*}

We remark that $\bm{\Sigma}_h^*$ is a vector field, we can extend it to define the stress space, i.e.,  $\underline{\Sigma}_h:=[\bm{\Sigma}_h^*]^2$.
We specify the degrees of freedom for $\bm{\Sigma}_h^*$ as follows: for $e\in \mathcal{F}_{pr}^{\frac{1}{2}}$
\begin{align}
	\phi_e(\bm{v}):=\int_e \bm{v}\cdot\bm{n}.\label{eq:dof}
\end{align}

\begin{lemma}
	$\bm{\Sigma}_h^*$ is uniquely determined by the degrees of freedom given in \eqref{eq:dof}.
	
\end{lemma}

\begin{proof}
	Since $\bm{v}$ is a vector function with two components and each component is a constant over each subcell, then we have
	\begin{align*}
		\text{dim}(\bm{\Sigma}_h^*)=2|\mathcal{T}_h|-|\mathcal{F}_{pr}^{\frac{1}{2}}\backslash \partial \Omega|.
	\end{align*}
	Let $\text{SD}$ represent the degrees of freedom corresponding to \eqref{eq:dof}, then we have
$
		|\text{SD}|=|\mathcal{F}_{pr}^{\frac{1}{2}}|.
$
	We can associate each edge in $\mathcal{F}_{pr}^{\frac{1}{2}}\backslash \partial \Omega$ to two subcells in $\mathcal{T}_h$ and associate each edge in $\partial \Omega$ to one subcell in $\mathcal{T}_h$, therefore, we have
$
		\text{dim}(\bm{\Sigma}_h^*)-|\text{SD}|=0.
$
	Then it suffices to show the uniqueness. Suppose $\bm{v}\in \bm{\Sigma}_h^*$ is defined such that all degrees of freedom corresponding to \eqref{eq:dof} are equal to zero. Then we have for $\bm{v}$ restricted to each subcell $E$ satisfying
$
		\bm{v}\cdot\bm{n}_1=0,\; \bm{v}\cdot\bm{n}_2=0,
$
	where $\bm{n}_i,i=1,2,$ are the unit normal vector of the edge in $\mathcal{F}_{pr}^{\frac{1}{2}}\cap \partial E$. Then we can infer that $\bm{v}=\bm{0}$. Thus, the proof is completed.
	
\end{proof}

To favor later analysis, we define the following semi-norm/norm
\begin{align*}
	\|\bm{v}_h\|_{1,h}^2= \sum_{e\in \mathcal{F}_{pr}^{\frac{1}{2}}}h_e^{-1}\|\jump{\bm{v}_h}\|_{L^2(e)}^2,\quad \|(\bm{v}_h,\eta_h)\|_{h}^2=\|\bm{v}_h\|_{1,h}^2+\|\eta_h\|_{L^2(\Omega)}^2.
\end{align*}

\subsection{The discrete formulation and statement of the main results}

In the following we use $\Gamma_h$ to represent $\Gamma_h^1$ and $\Gamma_h^2$ if there is no need to distinguish them. 
The discrete formulation for \eqref{eq:elasticity1}-\eqref{eq:elasticity2} reads as follows: Find $(\underline{\sigma}_h,\bm{u}_h, \gamma_h)\in \underline{\Sigma}_h\times \bm{U}_h\times \Gamma_h$ such that
\begin{alignat}{2}
	(\mathcal{A}\underline{\sigma}_h,\underline{w}_h)-\sum_{e\in \mathcal{F}_{dl}} (\bm{u}_h,\jump{\underline{w}_h\bm{n}})_e+(\text{as}(\underline{w}_h), \gamma_h)&=0 &&\quad \forall \underline{w}_h\in \underline{\Sigma}_h,\label{eq:MPFA1}\\
	-\sum_{e\in \mathcal{F}_{pr}^{\frac{1}{2}}}(\underline{\sigma}_h\bm{n},\jump{\bm{v}_h})_e&=(\bm{f},\bm{v}_h)&&\quad \forall \bm{v}_h\in \bm{U}_h\label{eq:MPFA2},\\
	(\text{as}(\underline{\sigma}_h), \xi_h)&=0&&\quad \forall \xi_h\in \Gamma_h. \label{eq:MPFA3}
\end{alignat}
We remark that the first term on the left-hand side of \eqref{eq:MPFA1} is block-diagonal, thereby, we can locally eliminate $\underline{\sigma}_h$ to get a reduced system. This is further explained in section~\ref{sec:cell-center}. 

\begin{remark}(local conservation).  For each macro-element $M$, setting $\bm{v}=\bm{1}$ in $M$ and zero otherwise, we can achieve the following local conservation owing to \eqref{eq:MPFA2}
	\begin{align*}
		-(\underline{\sigma}_h\bm{n},\bm{1})_{\partial M}&=\int_M \bm{f}.
	\end{align*}
	
This local conservation is important in many applications. Moreover, it is also crucial for constructing equilibrated stress in the context of a posteriori error estimation.
	\label{re:local}
\end{remark}

\begin{remark}(link to existing methods). Our method is closely related to the method proposed in \cite{Ambartsumyan20simplicial,Ambartsumyan21}. In \cite{Ambartsumyan20simplicial,Ambartsumyan21}, $\textnormal{BDM}_1$ is used for the approximation of stress, and piecewise constant function is used for the approximation of displacement, then a vertex quadrature rule is invoked to localize the stress degrees of freedom around each vertex, which allows the local elimination of the stress, resulting in a cell-centered system.  Our methods use piecewise constant functions for all the invoked unknowns. Instead of using the vertex quadrature rule, we choose the stress space to be broken $H(\tdiv;\Omega)$-conforming, which can naturally localize the stress degrees of freedom. As such, the stress bilinear form is block-diagonal, which resembles the methods proposed in \cite{Ambartsumyan20simplicial,Ambartsumyan21}.
	
	On the other hand, we also compare our method with the staggered cell-centered DG method developed by L. Zhao and E.-J. Park \cite{Zhao20}. In \cite{Zhao20}, the method is constructed by connecting the interior points to the vertices of the macro-element, and the stress space is defined to be normal continuous over the dual edges (our current method is normal continuous over the primal half edges), which is the key difference to our current setting. Both the method proposed in \cite{Zhao20} and our current method allow local elimination of the stress.
	
\end{remark}

\begin{remark}
	
	For problems with discontinuous
	compliance tensor $\mathcal{A}$, the rotation maybe discontinuous across the elements sharing the vertex, then 
	the rotation space in Method 1 may result in huge error near the interface of high and low 
	coefficient  regions. To this end, 
	we consider a modified variational formulation based on the scaled rotation $\tilde{\gamma}=\mathcal{A}^{-1}\gamma$ by 
	observing $\underline{\sigma}=\mathcal{A}^{-1}\nabla \bm{u}-\mathcal{A}^{-1}\gamma$. 
	 If $\mathcal{A}$ is strong discontinuous, $\underline{\sigma}$ is smoother than $\mathcal{A}\underline{\sigma}$ which means $\tilde{\gamma}$ is 
	 smoother than $\gamma$.
	 As a result, the modified weak formulation is to find $(\underline{\sigma}_h,\bm{u}_h, \tilde{\gamma}_h)\in \underline{\Sigma}_h\times \bm{U}_h\times \Gamma_h$ such that
	\begin{alignat}{2}
		(\mathcal{A}\underline{\sigma}_h,\underline{w}_h)-\sum_{e\in \mathcal{F}_{dl}} (\bm{u}_h,\jump{\underline{w}_h\bm{n}})_e+
		(\textnormal{as}(\mathcal{A}\underline{w}_h), \tilde{\gamma}_h)&=0 &&\quad \forall \underline{w}_h\in \underline{\Sigma}_h,\label{eq:MPFA1b}\\
		-\sum_{e\in \mathcal{F}_{pr}^{\frac{1}{2}}}(\underline{\sigma}_h\bm{n},\jump{\bm{v}_h})_e&=(\bm{f},\bm{v}_h)&&\quad \forall \bm{v}_h\in \bm{U}_h\label{eq:MPFA2b},\\
		(\textnormal{as}(\mathcal{A}\underline{\sigma_h}), \xi_h)&=0&&\quad \forall \xi_h\in \Gamma_h. \label{eq:MPFA3b}
	\end{alignat}
	By the definition of $\mathcal{A}$, $\textnormal{as}(\mathcal{A}\underline{\sigma_h})=1/(2\mu)\textnormal{as}(\underline{\sigma}_h)$.
	We will provide an example with highly heterogeneous $\mathcal{A}$  to illustrate the performance of \eqref{eq:MPFA1b}-\eqref{eq:MPFA3b}.
	To have a uniform presentation in relation to Method 2, in the following, we carry out 
	the well-posedness and error analysis for \eqref{eq:MPFA1}-\eqref{eq:MPFA3}. The
	analysis for the modified formulation \eqref{eq:MPFA1b}-\eqref{eq:MPFA3b} is similar.
\end{remark}


Before closing this section, we state the main results of the paper. 

		Let $(\underline{\sigma},\bm{u},\gamma)$ be the exact solution of \eqref{eq:elasticity1}-\eqref{eq:elasticity2} and let $(\underline{\sigma}_h,\bm{u}_h,\gamma_h)\in \underline{\Sigma}_h\times \bm{U}_h\times \Gamma_h$ be the discrete solution of \eqref{eq:MPFA1}-\eqref{eq:MPFA3} obtained by either Method 1 or Method 2. Assume that $(\underline{\sigma},\bm{u})\in \underline{H}^1(\Omega)\times \bm{H}^1(\Omega)$, then there exists a positive constant $C$ independent of the mesh size such that
			\begin{align*}
	\|\underline{\sigma}-\underline{\sigma}_h\|_{L^2(\Omega)}&\leq C  h\Big( \|\underline{\sigma}\|_{H^1(\Omega)}+\|\gamma\|_{H^1(\Omega)}\Big),\\
		\|\bm{u}-\bm{u}_h\|_{L^2(\Omega)}+\|\gamma-\gamma_h\|_{L^2(\Omega)}&\leq Ch\Big( \|\underline{\sigma}\|_{H^1(\Omega)}+\|\gamma\|_{H^1(\Omega)}\Big).
\end{align*}
In addition, if we assume  $\tdiv \underline{\sigma}\in \bm{H}^1(\Omega)$, then the following superconvergence holds
				\begin{align*}
	\|\bm{Q}_h\bm{u}-\bm{u}_h\|_{L^2(\Omega)}\leq C h^2\Big(\|\underline{\sigma}\|_{H^1(\Omega)}+\|\gamma\|_{H^1(\Omega)}+\|\tdiv\underline{\sigma}\|_{H^1(\Omega)}\Big),
\end{align*}
where $\bm{Q}_h$ is defined in \eqref{eq:Qh}.

\section{Error analysis}\label{sec:analysis}
In this section, we show the unique solvability of the discrete formulation (cf. Theorem~\ref{thm:unique}) and prove the convergence error estimates for $L^2$-error of stress and rotation; see Theorem~\ref{thm:uL2}. 

\subsection{Some elementary properties}

To begin, we introduce the bilinear bijection mapping to facilitate later analysis. For any element $E\in \mathcal{T}_h$ there exists a  bilinear bijection mapping $F: \widehat{E}\rightarrow E$, where $\widehat{E}=[0,1]^2$ is the reference square. Denote the Jacobian matrix by $\text{DF}$ and let $J=|\text{det}(\text{DF})|$.  For $\bm{x}=F(\widehat{x})$, we have
\begin{align*}
	\text{DF}^{-1}(\bm{x})=(\text{DF})^{-1}(\widehat{\bm{x}}),\quad J_{F^{-1}}=\frac{1}{J(\widehat{\bm{x}})}.
\end{align*}
Let $\widehat{E}$ has vertices $\widehat{r}_1=(0,0)^T,\widehat{r}_2=(1,0)^T, \widehat{r}_3=(1,1)^T$ and $\widehat{r}_4=(0,1)^T $. The bilinear mapping $F$ and its Jacobian matrix are given by
\begin{align*}
	F(\widehat{r})&=\bm{r}_1+\bm{r}_{21}\widehat{x}+\bm{r}_{41}\widehat{y}+(\bm{r}_{34}-\bm{r}_{21})\widehat{x}\widehat{y},\\\text{DF} &=[\bm{r}_{21},\bm{r}_{41}]+[(\bm{r}_{34}-\bm{r}_{21}) \widehat{y}, (\bm{r}_{34}-\bm{r}_{21})\widehat{x}],
\end{align*}
where $\bm{r}_{ij}=\bm{r}_i-\bm{r}_j$. In addition, we recall the Piola transformation defined as follows
\begin{align*}
	\bm{v}=\frac{1}{J}\text{DF} \widehat{\bm{v}}\circ F^{-1}.
\end{align*}
$\|\bm{v}\|_{L^2(M)}$ and $\|\widehat{\bm{v}}\|_{L^2(M)}$ are equivalent uniformly in $h$. Moreover, it is shown in \cite[Lemma~5.5]{Ewing99} that
\begin{align*}
	|\widehat{\bm{v}}|_{j,\widehat{M}}\leq C h^j \|\bm{v}\|_{j,M}, \quad j\geq 0.
\end{align*}

 For the rest of the paper, we assume that the quadrilateral elements are $\mathcal{O}(h^2)$-perturbations of parallelograms, following terminology from \cite{Ewing99}. Elements of this type are obtained by uniform refinements of  general quadrilateral grids. In this case, we have
\begin{align*}
	|\text{DF}|_{W^{1,\infty}(M)}\leq C h^2\quad \mbox{and}\quad |\frac{1}{J} \text{DF}|_{W^{j,\infty}(M)}\leq C h^{j-1},\quad j=1,2. 
\end{align*}

We let $P_h^{i}$ be the $L^2$-orthogonal projection onto $\Gamma_h,i=1,2$ corresponding to Method 1 and Method 2. More specifically, for Method 1, we have
\begin{align*}
	(P_h q, \mu)_D = (q,\mu)_D\quad \forall \mu\in P_0(D), D\in \mathcal{T}_D
\end{align*}
and for Method 2, we have
\begin{align*}
	(P_h q, \mu)_M= (q,\mu)_M \quad \forall \mu\in P_0(M), M\in \mathcal{T}_M.
\end{align*}
The following approximation properties holds for $P_h^i$
\begin{align*}
	\|q-P_h^1q\|_{L^2(D)}\leq C h_D\|\nabla q\|_{L^2(D)}\quad \forall q\in H^1(D),
	\|q-P_h^2q\|_{L^2(M)}\leq C h_M\|\nabla q\|_{L^2(M)}\quad \forall q\in H^1(M).
\end{align*}
We use $\bm{Q}_h$ to represent the $L^2$-orthogonal projection onto $\bm{U}_h$, i.e., 
\begin{align}
	(\bm{Q}_h\bm{v}, \bm{\theta})_M = (\bm{v}, \bm{\theta})_M\quad \forall \bm{\theta}\in \bm{P}_0(M),M\in \mathcal{T}_{M}.\label{eq:Qh}
\end{align}
It holds
\begin{align*}
	\|\bm{v}-\bm{Q}_h\bm{v}\|_{L^2(M)}\leq C h_M \|\bm{v}\|_{H^1(M)}\quad \forall \bm{v}\in \bm{H}^1(M), M\in \mathcal{T}_M.
\end{align*}
We define for $E\in \mathcal{T}_h$ and $\bm{v}\in \bm{H}^1(E)$
\begin{align}
	(\bm{\mathcal{J}}\bm{v}\cdot\bm{n}, \mu)_e=(\bm{v}\cdot\bm{n},\mu)_e\quad \forall \mu\in P_0(e),e\in \mathcal{F}_{pr}^{\frac{1}{2}}\cap \partial E,\label{eq:Jhdef}
\end{align}
which is well-defined owing to the degrees of freedom given in \eqref{eq:dof}. It is easy to check that $\bm{\Pi}^u\bm{v}=\bm{v}$ if $\bm{v}$ is a constant vector over each macro-element $M$, then an application of the Bramble-Hilbert lemma yields
\begin{align*}
	\|\bm{v}-\bm{\mathcal{J}}\bm{v}\|_{L^2(M)}\leq C h |\bm{v}|_{H^1(M)}\quad \forall \bm{v}\in \bm{H}^1(M), M\in \mathcal{T}_M.
\end{align*}
We denote the lowest order Raviart-Thomas space on rectangular mesh as $\bm{U}_h^{\text{RT}}$ (cf. \cite{Brezzi91}). 
The projection operator $\bm{\Pi}_h:\bm{H}^1(M)\rightarrow \bm{U}_h^{\text{RT}}$ is defined such that the following holds
\begin{align}
	(\bm{\Pi}_h\bm{v}\cdot\bm{n},q)_e=(\bm{v}\cdot\bm{n},q)_e\quad \forall q\in P_0(e),\bm{v}\in \bm{H}^1(M), M\in \mathcal{T}_M, s\subset \partial M.\label{eq:Pihdef}
\end{align}
When $\bm{v}\in \bm{U}_h$, this means $\bm{v}\notin \bm{H}^1(M)$, \eqref{eq:Pihdef} can be interpreted as 
\begin{align*}
	(\bm{\Pi}_h\bm{v}\cdot\bm{n},q)_e=\sum_{k=1}^2(\bm{v}\cdot\bm{n},q)_{e_k},
\end{align*}
where $e_k\in \mathcal{F}_{pr}^{\frac{1}{2}}$ is the half edge of $e$. In addition, the following holds (cf. \cite[(3.9)]{Wheeler06})
\begin{align*}
	\|\bm{\Pi}_h\bm{v}\|_{H^1(M)}\leq C \|\bm{v}\|_{H^1(M)}\quad \forall \bm{v}\in \bm{H}^1(M),M\in \mathcal{T}_M.
\end{align*}
Furthermore, the following error estimates also hold
\begin{align}
	\|\bm{v}-\bm{\Pi}_h\bm{v}\|_{L^2(M)}&\leq C h_M |\bm{v}|_{H^1(M)}\quad \forall \bm{v}\in H^1(M),M\in \mathcal{T}_M,\label{eq:Ph1}\\
	\|\nabla\cdot (\bm{v}-\bm{\Pi}_h\bm{v})\|_{L^2(M)}&\leq Ch^r \|\nabla\cdot\bm{v}\|_{H^r(M)},0\leq r\leq 1.\label{eq:Ph2}
\end{align}
\Red{For a vector field $\bm{v}:=(v_1,v_2)^T$, we define
		\begin{align*}
		\tcurl\bm{v}=\begin{pmatrix}
			\partial_y v_1 & -\partial_xv_1\\
			\partial_yv_2 & -\partial_xv_2
		\end{pmatrix}.
	\end{align*}
}

We use $\underline{I}_h^e$ and $\underline{\Pi}_h$ to represent the extension of $\bm{\mathcal{J}}$ and $\bm{\Pi}_h$ to the matrix field, i.e., $\underline{I}_h^e=[\bm{\Pi}^u]^2$ and $\underline{\Pi}_h=[\bm{\Pi}_h]^2$. We can deduce from the definition of $\underline{\Pi}_h$ that
\begin{align*}
	(\bm{v}_h,\tdiv\underline{\Pi}_h\underline{w}_h)=\sum_{M\in \mathcal{T}_M}(\bm{v}_h, \underline{\Pi}_h\underline{w}_h \bm{n})_{\partial M}=\sum_{M\in \mathcal{T}_M}(\bm{v}_h,\underline{w}_h\bm{n})_{\partial M}=\sum_{e\in \mathcal{F}_{pr}^{\frac{1}{2}}} (\jump{\bm{v}_h},\underline{w}_h\bm{n})_e\quad \forall \bm{v}_h\in \bm{U}_h,\underline{w}_h\in \underline{\Sigma}_h.
\end{align*}
Therefore, the discrete formulation can be recast into the following equivalent form: Find $(\underline{\sigma}_h,\bm{u}_h,\gamma_h)\in \underline{\Sigma}_h\times \bm{U}_h\times \Gamma_h$ such that
\begin{alignat}{2}
	(\mathcal{A}\underline{\sigma}_h,\underline{w}_h)+(\bm{u}_h,\tdiv\underline{\Pi}_h\underline{w}_h)+(\text{as}(\underline{w}_h), \gamma_h)&=0&&\quad \forall \underline{w}_h\in \underline{\Sigma}_h,\label{eq:MPFA1-RT}\\
	-(\tdiv\underline{\Pi}_h \underline{\sigma}_h,\bm{v}_h)&=(\bm{f},\bm{v}_h)&&\quad \forall \bm{v}_h\in \bm{U}_h\label{eq:MPFA2-RT},\\
	(\text{as}(\underline{\sigma}_h), \xi_h)&=0&&\quad \forall \xi_h\in \Gamma_h.
\end{alignat}
It follows from \eqref{eq:MPFA2-RT} that
\begin{align}
	\tdiv\underline{\Pi}_h \underline{\sigma}_h=-\bm{Q}_h\bm{f}.\label{eq:mass-sigma}
\end{align}
We can observe from \eqref{eq:mass-sigma} that we can construct the $\underline{H}(\tdiv;\Omega)$-conforming stress by simply applying the projection operator $\underline{\Pi}_h$, which is computationally cheap.

%
%
%
%
%
%

Now we recall the following stabilized $\bm{Q}_1-Q_0$ pair for the Stokes equations (cf. \cite{Kechkar92}).

\begin{lemma}\label{lemma:stokes}
	We define the finite element spaces for each macro-element $M,M\in \mathcal{T}_{pr}$ as follows:
	\begin{align*}
		\bm{V}_h(M):&=\{\bm{q}=(q_1,q_2)^T\in \bm{H}_0^1(M): \bm{q}_{|E}\in \bm{Q}_1(E),\forall E\in \mathcal{T}_h\; \textnormal{such that} \; E\subset M\},\\
		W_h(M):&=\{w\in L^2(\Omega): w_{|E} \in Q_0(E),\forall E\in \mathcal{T}_h\; \textnormal{such that} \;E \subset M\},
	\end{align*}
Then, there exists a unique solution $(\bm{\rho}_h, p_h^*)\in \bm{V}_h(M)\times W_h(M)$ to the following system
\begin{align}
	(\nabla \bm{\rho}_h,\nabla \bm{v})_M+(p_h^*,\nabla\cdot \bm{v})_M&=(\bm{f},\bm{v})_M\quad \forall \bm{v}\in \bm{V}_h(M),\label{eq:stokes1}\\
	(\tdiv \bm{\rho}_h,q)_M+\beta\sum_{e\in \mathcal{F}_{dl}\cap M}h_e(\jump{p_h^*},\jump{q})_e&=(g,q)_M\quad \forall q\in W_h(M)\label{eq:stokes2},
\end{align}
where $\beta>0$ is the stabilization parameter. 
Moreover, it holds
\begin{align}
	\|\bm{\rho}_h\|_{H^1(\Omega)}\leq C \|g\|_{L^2(\Omega)}\label{eq:stokes-stability},
\end{align}
where $\bm{\rho}_h$ over each macro-element is defined by \eqref{eq:stokes1}-\eqref{eq:stokes2}. 

\end{lemma}

\begin{lemma}(inf-sup condition).\label{eq:inf-sup-elasticity} For Method 2, 
	there exists a positive constant $C$ independent of the mesh size such that
		\begin{align}
		\sup_{\underline{0}\neq \underline{w}_h\in \underline{\Sigma}_h}	\frac{\sum_{e\in \mathcal{F}_{pr}^{\frac{1}{2}}}(\underline{w}_h\bm{n},\jump{\bm{v}_h})_e+(\textnormal{as}(\underline{w}_h),\eta_h)}{\|\underline{w}_h\|_{L^2(\Omega)}}\geq C \|(\bm{v}_h,\eta_h)\|_{h}\quad \forall (\bm{v}_h, \eta_h)\in \bm{U}_h\times \Gamma_h.	\label{eq:inf}
	\end{align}

\end{lemma}

\begin{proof}

	Over each subcell $E\in \mathcal{T}_h$, we define 
		\begin{align}
		(\underline{\psi}_1\bm{n}, \bm{\mu})_e=h_e^{-1}(\jump{\bm{v}_h}, \bm{\mu})_e \quad \forall \bm{\mu}\in\bm{P}_0(e),e\in \mathcal{F}_{pr}^{\frac{1}{2}}\cap \partial E,\label{eq:vhdef}
	\end{align}
	which implies that
	\begin{align*}
		\sum_{e\in \mathcal{F}_{pr}^{\frac{1}{2}}}(\underline{\psi}_1\bm{n},\jump{\bm{v}_h})_e=\sum_{e\in \mathcal{F}_{pr}^{\frac{1}{2}}}h_e^{-1}\|\jump{\bm{v}_h}\|_{L^2(e)}^2=\|\bm{v}_h\|_{1,h}^2.
	\end{align*}
	Moreover, the scaling arguments and \eqref{eq:vhdef} imply
	\begin{align}
		\|\underline{\psi}_1\|_{L^2(\Omega)}\leq C \|\bm{v}_h\|_{1,h}.\label{eq:psi1vh}
	\end{align}
We define $\underline{\psi}_2\in \underline{\Sigma}_h$ by
	\begin{align}
\Red{\underline{\psi}_2=\tcurl(\bm{\rho}_h)\label{eq:psi2n}.}
	\end{align}
Since $(\bm{\rho}_h)_{|M}\in \bm{V}_h(M)$, it is easy to check that $\tcurl (\bm{\rho}_h)\in \underline{\Sigma}_h$.  
	Therefore, in view of \eqref{eq:stokes2} with $g=\eta_h-\text{as}(\underline{\psi}_1)$, we can obtain
\begin{align*}
(\text{as}(\underline{\psi}_1+\underline{\psi}_2),\eta_h)=(\text{as}(\underline{\psi}_1),\eta_h)+(\tdiv \bm{\rho}_h, \eta_h)=\|\eta_h\|_{L^2(\Omega)}^2.
\end{align*}
Here we use the fact that $\eta_h$ is a constant function over each macro-element $M$, and thus, the jump term for pressure vanishes. 

It follows from \eqref{eq:stokes-stability}  and \eqref{eq:psi1vh} that
\begin{align}
	\|\underline{\psi}_2\|_{L^2(\Omega)}\leq C\|\bm{\rho}_h\|_{H^1(\Omega)}\leq C \|\eta_h-\text{as}(\underline{\psi}_1)\|_{L^2(\Omega)}\leq C \|(\bm{v}_h,\eta_h)\|_{h}.\label{eq:psi2bound}
\end{align}
Then, we define $\underline{w}_h=\underline{\psi}_1+\underline{\psi}_2$, it will lead to
\begin{align*}
	\sum_{e\in \mathcal{F}_{pr}^{\frac{1}{2}}}(\underline{w}_h\bm{n},\jump{\bm{v}_h})_e+(\text{as}(\underline{w}_h),\eta_h)
	=\sum_{e\in \mathcal{F}_{pr}^{\frac{1}{2}}}h_e^{-1}\|\jump{\bm{v}_h}\|_{L^2(e)}^2+\|\eta_h\|_{L^2(\Omega)}^2,
\end{align*}
where we use $\tcurl (\bm{\rho}_h)\bm{n}_{|e}=\bm{0}$ for any $e\in \mathcal{F}_{pr}^{\frac{1}{2}}$ if $\bm{\rho}_h=\bm{0}$ on $\partial M$.
This, \eqref{eq:psi1vh} and \eqref{eq:psi2bound} complete the proof. 
	
\end{proof}

	\begin{lemma}(inf-sup condition).\label{lemma:inf-sup2}
	For Method 2, there exists a positive constant $C$ such that
		\begin{align*}
			\sup_{\bm{0}\neq \underline{w}\in \underline{\Sigma}_h}\frac{(\bm{v},\tdiv \underline{\Pi}_h\underline{w})+(\textnormal{as}(\underline{w}),\eta)}{\|\underline{w}\|_{L^2(\Omega)}}\geq C (\|\bm{v}\|_{L^2(\Omega)}+\|\eta\|_{L^2(\Omega)})\quad \forall (\bm{v},\eta)\in \bm{U}_h\times \Gamma_h.
		\end{align*}

	\end{lemma}
	
	\begin{proof}
		For $\bm{v}\in \bm{L}^2(\Omega)$, there exists $\underline{\theta}\in \underline{H}^1(\Omega)$ (see, for example, \cite{Girault79}), such that 
		\begin{align}
			\tdiv\underline{\theta}=\bm{v},\quad
			\|\underline{\theta}\|_{H^1(\Omega)}\leq C \|\bm{v}\|_{L^2(\Omega)}.\label{eq:divtheta}
		\end{align}
		It follows from the trace inequality and the interpolation error estimate that
		\begin{align}
			\|\underline{I}_h^e\underline{\theta}\|_{L^2(\Omega)}\leq C \|\underline{\theta}\|_{H^1(\Omega)}.\label{eq:Jhtheta1h}
		\end{align}
		We can infer from \eqref{eq:divtheta} and \eqref{eq:Jhdef} that
		\begin{align*}
			\|\bm{v}\|_{L^2(\Omega)}^2=(\bm{v},\tdiv \underline{\theta})=\sum_{e\in \mathcal{F}_{pr}^{\frac{1}{2}}} (\jump{\bm{v}}, \underline{\theta}\bm{n})_e=\sum_{e\in \mathcal{F}_{pr}^{\frac{1}{2}}} (\jump{\bm{v}},\underline{I}_h^e \underline{\theta}\bm{n})_e.
		\end{align*}
		Then  we have
		\begin{equation}
			\begin{split}
				(\bm{v}, \tdiv\underline{\Pi}_h \underline{I}_h^e\underline{\theta})&=\sum_{M\in \mathcal{T}_M} (\bm{v}, \underline{\Pi}_h\underline{I}_h^e\underline{\theta}\bm{n})_{\partial M}\\
				&=\sum_{M\in \mathcal{T}_M} (\bm{v}, \underline{I}_h^e\underline{\theta}\bm{n})_{\partial M}=\sum_{e\in \mathcal{F}_{pr}^{\frac{1}{2}}} (\jump{\bm{v}}, \underline{I}_h^e\underline{\theta}\bm{n})_e=\|\bm{v}\|_{L^2(\Omega)}^2.
			\end{split}
			\label{eq:qJhtheta}
		\end{equation}
		We define $\underline{w}_2=\tcurl(\bm{\rho}_h)$ with $g=\eta-\underline{I}_h^e\underline{\theta}$ (cf. Lemma~\ref{lemma:stokes}). Then  we have $\tdiv \underline{w}_2=\bm{0}$ and an application of Lemma~\ref{lemma:stokes} yields
		\begin{align*}
			(\text{as}(\underline{I}_h^e\underline{\theta}+\underline{w}_2),\eta)=\|\eta\|_{L^2(\Omega)}^2\quad \mbox{and}\quad \|\underline{w}_2\|_{L^2(\Omega)}\leq C \|\bm{v}\|_{L^2(\Omega)}.
		\end{align*}
		This and \eqref{eq:qJhtheta} complete the proof by taking $\underline{w}=\underline{I}_h^e\underline{\theta}+\underline{w}_2$.

	\end{proof}

\begin{remark}\label{remark:inf-sup}
	We remark that the proof for Lemma~\ref{eq:inf-sup-elasticity} and Lemma~\ref{lemma:inf-sup2} can be easily extended to proving the inf-sup condition for Method 1. The only difference is to define the Stokes pair used in Lemma~\ref{lemma:stokes} for each $D,D\in \mathcal{T}_D$; indeed, we treat $D$ as the ``macro-element" in this case. 
	
\end{remark}

%
%
%
%
%
%
%
%

\begin{theorem}\label{thm:unique}(existence and uniqueness).
	There exists a unique solution to \eqref{eq:MPFA1}-\eqref{eq:MPFA2}. In addition, there exists a positive constant $C$ independent of the mesh size such that the following estimate holds
	\begin{align}
		\|\underline{\sigma}_h\|_{L^2(\Omega)}+\|(\bm{u}_h,\gamma_h)\|_{h}\leq C\|\bm{f}\|_{L^2(\Omega)}.\label{eq:stability}
	\end{align}

\end{theorem}

\begin{proof}
	
	We first show the stability estimate \eqref{eq:stability}. 
	Taking $\underline{w}_h=\underline{\sigma}_h$, $\bm{v}_h=\bm{u}_h,\xi_h= \gamma_h$ in  \eqref{eq:MPFA1}-\eqref{eq:MPFA2} and summing up the resulting equations, we have
	\begin{align}
		\|\underline{\sigma}_h\|_{\mathcal{A}}^2=(\bm{f}, \bm{u}_h).\label{eq:sigmaA}
	\end{align}
Next, it follows from Lemma~\ref{eq:inf-sup-elasticity} and \eqref{eq:MPFA1} that
	\begin{align*}
		\|(\bm{u}_h,\gamma_h)\|_{h}\leq C 	\|\mathcal{A}\underline{\sigma}_h\|_{L^2(\Omega)}.
	\end{align*}
This,  \eqref{eq:sigmaA}, the Cauchy-Schwarz inequality and the discrete Poincar\'{e} inequality give  (cf. \cite{Brenner03}) 
\begin{align*}
		\|\underline{\sigma}_h\|_{\mathcal{A}}&\leq C \|\bm{f}\|_{L^2(\Omega)} \|\bm{u}_h\|_{L^2(\Omega)}\leq C \|\bm{f}\|_{L^2(\Omega)} \|\bm{u}_h\|_{1,h}\leq C\|\bm{f}\|_{L^2(\Omega)} 	\|\underline{\sigma}_h\|_{\mathcal{A}},
\end{align*}
which combined with Young's inequality implies
\begin{align*}
	\|\underline{\sigma}_h\|_{\mathcal{A}}\leq C\|\bm{f}\|_{L^2(\Omega)}.
\end{align*}
As such, we have $\|(\bm{u}_h,\gamma_h)\|_{h}\leq C\|\bm{f}\|_{L^2(\Omega)}$.


Since \eqref{eq:MPFA1}-\eqref{eq:MPFA2} is a square linear system, the uniqueness of solution implies the existence. Setting $\bm{f}=\bm{0}$ in \eqref{eq:stability} implies $\underline{\sigma}_h=\underline{0},  \Red{\bm{u}_h=\bm{0}}$ and $\gamma_h=0$.

%

\end{proof}

The trapezoidal quadrature rule is exact for scalar linear functions, therefore, we have
\begin{align*}
	\int_{\hat{E}} \hat{v}_k = \frac{1}{4}\sum_{i=1}^4 \hat{v}_k\mid e_{ik}, \quad k=1,2,
\end{align*}
where $e_{ik}$ is the half edge of $e\in \partial M$.

				\begin{lemma}\label{lemma:average}
					
					For $\bm{\chi}\in \bm{P}_0(\widehat{M})$, $M\in \mathcal{T}_M$, where $\widehat{M}$ is the reference element corresponding to $M$, it holds
					\begin{align*}
						(\bm{\chi}, \bm{v}-\bm{\Pi}_h\bm{v})_{\widehat{M}}=0 \quad \forall \bm{v}\in \bm{\Sigma}_h^*.
					\end{align*}
					
				\end{lemma}
				
				\begin{proof}
					Let $\bm{v}=(v_1,v_2)$, then the exactness of the trapezoidal quadrature rule for scalar linear functions imply that
					\begin{align*}
						\int_{\widehat{M}} (\bm{\Pi}_h\bm{v})_k=\frac{1}{4} \sum_{i=1}^4(\bm{\Pi}_h\bm{v})_k\mid_{e_{ik}}\quad k=1,2.
					\end{align*}
			As s consequence, it holds		
					\begin{align*}
						\int_{\widehat{M}} ((I-\bm{\Pi}_h)\bm{v})_k=0,
					\end{align*}
					which completes the proof.
					
				\end{proof}
			
				\begin{lemma}\label{lemma:phi}
				For any $\underline{w}_h\in \underline{\Sigma}_h$ and $\underline{\theta}\in \underline{H}^1(\Omega)$, we have
				\begin{align*}	(\underline{\Pi}_h\underline{w}_h-\underline{w}_h,\underline{\theta})\leq C h \|\underline{\Pi}_h\underline{w}_h-\underline{w}_h\|_{L^2(\Omega)} |\underline{\theta}|_{H^1(\Omega)}.
				\end{align*}

			\end{lemma}
			
			\begin{proof}
				
				On the reference element, we have
				\begin{align*}
					(\widehat{\underline{\Pi}}_h\widehat{\underline{w}}_h-\widehat{\underline{w}}_h,\widehat{\underline{\theta}})_{\widehat{M}}\leq C \|\widehat{\underline{\Pi}}_h\widehat{\underline{w}}_h-\widehat{\underline{w}}_h\|_{L^2(\widehat{M})} \|\widehat{\underline{\theta}}\|_{H^1(\widehat{M})}.
				\end{align*}
				Owing to Lemma~\ref{lemma:average}, we have
				\begin{align*}
					(\widehat{\underline{\Pi}}_h\widehat{\underline{w}}_h-\widehat{\underline{w}}_h,\widehat{\underline{\theta}})_{\widehat{M}}\leq C \|\widehat{\underline{\Pi}}_h\widehat{\underline{w}}_h-\widehat{\underline{w}}_h\|_{L^2(\widehat{M})} \|\nabla \widehat{\underline{\theta}}\|_{L^2(\widehat{M})}.
				\end{align*}
				Then the transformation back to the element $M$ yields the desired estimate.

			\end{proof}

%
%
%
%
%

%
%

	\subsection{The convergence error analysis}

		\begin{lemma}
	Let $(\bm{u},p)$ be the exact solution and let $(\bm{u}_h,p_h)\in \bm{U}_h\times P_h$ be the discrete solution of \eqref{eq:MPFA1-RT}-\eqref{eq:MPFA2-RT}. 
	Then,  the following error equations hold
	\begin{align}
		(\mathcal{A}(\underline{\sigma}-\underline{\sigma}_h),\underline{w})+(\bm{u}-\bm{u}_h,\tdiv  \underline{\Pi}_h \underline{w})+(\gamma-\gamma_h, \textnormal{as}(\underline{w}))&=(\mathcal{A}\underline{\sigma}+\gamma\underline{\delta},\underline{w}-\underline{\Pi}_h\underline{w})\quad \forall \underline{w}\in \underline{\Sigma}_h,\label{eq:error1-RT1}\\
		(\tdiv \underline{\Pi}_h(\underline{\sigma}-\underline{\sigma}_h),\bm{v})&=0\quad \forall \bm{v}\in P_h\label{eq:error2-RT2},\\
		(\textnormal{as}(\underline{\sigma}-\underline{\sigma}_h), \xi)&=0\quad \forall \xi\in \Gamma_h\label{eq:error3-RT}.
	\end{align}

\end{lemma}

\begin{proof}
	First, we have from integration by parts that
	\begin{align*}
		(\bm{u},\tdiv \underline{\Pi}_h\underline{w})=-(\nabla \bm{u}, \underline{\Pi}_h\underline{w})=-(\mathcal{A}\underline{\sigma}+\gamma\underline{\delta},\underline{\Pi}_h\underline{w}),
	\end{align*}
	thereby,
	\begin{align*}
		(\mathcal{A}\underline{\sigma},\underline{w})+(\bm{u},\tdiv \underline{\Pi}_h\underline{w})+(\gamma, \text{as}(w))=(\mathcal{A}\underline{\sigma}+\gamma\underline{\delta},\underline{w}-\underline{\Pi}_h\underline{w}),
	\end{align*}
	which combined with \eqref{eq:MPFA1-RT} implies \eqref{eq:error1-RT1}.
	
	Then \eqref{eq:error2-RT2} follows directly from the definition of $\underline{\Pi}_h$ and \eqref{eq:MPFA2-RT}.
\end{proof}

			\begin{theorem}\label{thm:uL2}
	Let $(\underline{\sigma},\bm{u})$ be the exact solution and let $(\underline{\sigma}_h,\bm{u}_h)\in \underline{\Sigma}_h\times \bm{U}_h$ be the discrete solution of \eqref{eq:MPFA1-RT}-\eqref{eq:MPFA2-RT}. Assume that $(\underline{\sigma},\gamma)\in \underline{H}^1(\Omega)\times H^1(\Omega)$, then the following convergence error estimates hold
	\begin{align*}
		\|\underline{\sigma}-\underline{\sigma}_h\|_{L^2(\Omega)}\leq C h\Big( \|\underline{\sigma}\|_{H^1(\Omega)}+\|\gamma\|_{H^1(\Omega)}\Big)
	\end{align*}
and
			\begin{align*}
		\|\bm{Q}_h\bm{u}-\bm{u}_h\|_{L^2(\Omega)}+\|\gamma_h-P_h\gamma\|_{L^2(\Omega)}\leq Ch\Big( \|\underline{\sigma}\|_{H^1(\Omega)}+\|\gamma\|_{H^1(\Omega)}\Big).
	\end{align*}

\end{theorem}

\begin{proof}
	
	Taking $\underline{w}=\underline{I}_h^e\underline{\sigma}-\underline{\sigma}_h$, $\bm{v}=\bm{Q}_h\bm{u}-\bm{u}_h$ and $\xi=P_h\gamma-\gamma_h$ in \eqref{eq:error1-RT1}-\eqref{eq:error3-RT}, then it holds
	\begin{align*}
		(\mathcal{A}(\underline{\sigma}-\underline{\sigma}_h),\underline{I}_h^e\underline{\sigma}-\underline{\sigma}_h)+(\bm{u}-\bm{u}_h,\tdiv\underline{\Pi}_h (\underline{I}_h^e\underline{\sigma}_h-\underline{\sigma}_h))&+(\gamma-\gamma_h, \text{as}(\underline{I}_h^e\underline{\sigma}-\underline{\sigma}_h))\\
		&=(\mathcal{A}\underline{\sigma}+\gamma\underline{\delta},\underline{I}_h^e\underline{\sigma}-\underline{\sigma}_h-\underline{\Pi}_h(\underline{I}_h^e\underline{\sigma}-\underline{\sigma}_h)),\\
		(\tdiv \underline{\Pi}_h (\underline{\sigma}-\underline{\sigma}_h),\bm{Q}_h\bm{u}-\bm{u}_h)&=0,\\
		(\text{as}(\underline{\sigma}-\underline{\sigma}_h), P_h\gamma-\gamma_h)&=0.
	\end{align*}
	The definition of $\underline{\Pi}_h$ (cf. \eqref{eq:Pihdef}) yields
	\begin{align*}
		(\tdiv \underline{\Pi}_h (\underline{\sigma}-\underline{\sigma}_h),\bm{Q}_h\bm{u}-\bm{u}_h)=	(\tdiv (\underline{\sigma}-\underline{\Pi}_h\underline{\sigma}_h),\bm{Q}_h\bm{u}-\bm{u}_h).
	\end{align*}
	We can infer from the definition of $\underline{\Pi}_h$ (cf. \eqref{eq:Pihdef}), integration by parts and \eqref{eq:Jhdef} that
	\begin{align*}
		(\bm{u}-\bm{u}_h,\tdiv\underline{\Pi}_h(\underline{I}_h^e\underline{\sigma}-\underline{\sigma}_h))&=(\bm{Q}_h\bm{u}-\bm{u}_h,\tdiv\underline{\Pi}_h (\underline{I}_h^e\underline{\sigma}-\underline{\sigma}_h))\\
		&=\sum_{T\in \mathcal{T}_M}(\bm{Q}_h\bm{u}-\bm{u}_h, \underline{\Pi}_h(\underline{I}_h^e\underline{\sigma}-\underline{\sigma}_h)\bm{n})_{\partial T}\\
		&=\sum_{T\in \mathcal{T}_M}(\bm{Q}_h\bm{u}-\bm{u}_h, (\underline{I}_h^e\underline{\sigma}-\underline{\Pi}_h\underline{u}_h)\bm{n})_{\partial T}\\
		&=\sum_{T\in \mathcal{T}_M}(\bm{Q}_h\bm{u}-\bm{u}_h, (\underline{\sigma}-\underline{\Pi}_h\underline{\sigma}_h)\bm{n})_{\partial T}\\
		&=(\bm{Q}_h\bm{u}-\bm{u}_h, \tdiv (\underline{\sigma}-\underline{\Pi}_h\underline{\sigma}_h)).
	\end{align*}
	In view of Lemma~\ref{lemma:average} and the Piola transformation, we have
	\begin{align*}
		(\mathcal{A}\underline{\sigma}+\gamma\underline{\delta},\underline{I}_h^e\underline{\sigma}-\underline{\sigma}_h-\underline{\Pi}_h(\underline{I}_h^e\underline{\sigma}-\underline{\sigma}_h))& \leq C h( \|\underline{\sigma}\|_{H^1(\Omega)}+\|\gamma\|_{H^1(\Omega)})\|\underline{I}_h^e\underline{\sigma}-\underline{\sigma}_h\|_{L^2(\Omega)}.
	\end{align*}
	Therefore, we have from the Cauchy-Schwarz inequality that
	\begin{equation}
		\begin{split}
		\|\underline{I}_h^e\underline{\sigma}-\underline{\sigma}_h\|_{L^2(\Omega)}^2&=(\underline{I}_h^e\underline{\sigma}-\underline{\sigma}, \underline{I}_h^e\underline{\sigma}-\underline{\sigma}_h)+
		(\mathcal{A}\underline{\sigma}+\gamma\underline{\delta},\underline{I}_h^e\underline{\sigma}-\underline{\sigma}_h-\underline{\Pi}_h(\underline{I}_h^e\underline{\sigma}-\underline{\sigma}_h))\\
		&-(\gamma-P_h\gamma,\text{as}(\underline{I}_h^e\underline{\sigma}-\underline{\sigma}_h) )-	(\text{as}(\underline{\sigma}-\underline{I}_h^e\underline{\sigma}), P_h\gamma-\gamma_h)\\
		&\leq C \Big(h( \|\underline{\sigma}\|_{H^1(\Omega)}+\|\gamma\|_{H^1(\Omega)})+\|P_h\gamma-\gamma_h\|_{L^2(\Omega)}\Big)\|\underline{I}_h^e\underline{\sigma}-\underline{\sigma}_h\|_{L^2(\Omega)}.
		\end{split}
	\label{eq:Ihsigma}
	\end{equation}
		We have from the inf-sup condition (cf. Lemma~\ref{lemma:inf-sup2}) and \eqref{eq:error1-RT1} that
	\begin{equation}
		\begin{split}
			\|\bm{Q}_h\bm{u}-\bm{u}_h\|_{L^2(\Omega)}+\|\gamma_h-P_h\gamma\|_{L^2(\Omega)}&\leq C \sup_{\underline{0}\neq \underline{w}\in \underline{\Sigma}_h} \frac{(\bm{Q}_h\bm{u}-\bm{u}_h, \tdiv \underline{\Pi}_h\underline{w})+(\text{as}(\underline{w}),P_h\gamma-\gamma_h)}{\|\underline{w}\|_{L^2(\Omega)}}\\
			&=C \sup_{\underline{0}\neq \underline{w}\in \underline{\Sigma}_h}\frac{-(\mathcal{A}	(\underline{\sigma}-\underline{\sigma}_h),\underline{w})+(\mathcal{A}\underline{\sigma}+\gamma\underline{\delta},\underline{w}-\underline{\Pi}_h\underline{w})-(\gamma-P_h\gamma, \text{as}(\underline{w}))}{\|\underline{w}\|_{L^2(\Omega)}}.
		\end{split}
		\label{eq:Qhp-ph}
	\end{equation}
	The Cauchy-Schwarz inequality and the discrete Poincar\'{e} inequality lead to
	\begin{align*}
		-(\mathcal{A}	(\underline{\sigma}-\underline{\sigma}_h),\underline{w})-(\gamma-P_h\gamma, \text{as}(\underline{w}))\leq C \Big(\|\underline{\sigma}-\underline{\sigma}_h\|_{\mathcal{A}}+\|\gamma-P_h\gamma\|_{L^2(\Omega)}\Big)\|\underline{w}\|_{L^2(\Omega)}.
	\end{align*}
	An appeal to Lemma~\ref{lemma:average} yields
	\begin{align*}
		(\mathcal{A}\underline{\sigma}+\gamma \underline{\delta},\underline{w}-\underline{\Pi}_h\underline{w})\leq C h (\|\underline{\sigma}\|_{H^1(\Omega)}+\|\gamma\|_{H^1(\Omega)})\|\underline{w}\|_{L^2(\Omega)}.
	\end{align*}
	Therefore, we have
	\begin{align}
		\|\bm{Q}_h\bm{u}-\bm{u}_h\|_{L^2(\Omega)}+\|\gamma_h-P_h\gamma\|_{L^2(\Omega)}&\leq C \Big(\|\underline{\sigma}-\underline{\sigma}_h\|_{\mathcal{A}}+\|\gamma-P_h\gamma\|_{L^2(\Omega)}+ h(\|\underline{\sigma}\|_{H^1(\Omega)}+\|\gamma\|_{H^1(\Omega)})\Big).\label{eq:Qhu}
	\end{align}
Combining the above with \eqref{eq:Ihsigma} implies
	\begin{align*}
		\|\underline{I}_h^e\underline{\sigma}-\underline{\sigma}_h\|_{L^2(\Omega)}\leq C h (\|\underline{\sigma}\|_{H^1(\Omega)}+\|\gamma\|_{H^1(\Omega)}).
	\end{align*}
	As a consequence of the triangle inequality and the interpolation error estimate, one has
	\begin{align*}
		\|\underline{\sigma}-\underline{\sigma}_h\|_{L^2(\Omega)}\leq \|\underline{\sigma}-\underline{I}_h^e\underline{\sigma}\|_{L^2(\Omega)}+\|\underline{I}_h^e\underline{\sigma}-\underline{\sigma}_h\|_{L^2(\Omega)}\leq C h(\|\underline{\sigma}\|_{H^1(\Omega)}+\|\gamma\|_{H^1(\Omega)}).
	\end{align*}
This and \eqref{eq:Qhu} give
		\begin{align*}
	\|\bm{Q}_h\bm{u}-\bm{u}_h\|_{L^2(\Omega)}+\|\gamma_h-P_h\gamma\|_{L^2(\Omega)}&\leq  C h(\|\underline{\sigma}\|_{H^1(\Omega)}+\|\gamma\|_{H^1(\Omega)}).
\end{align*}
	Therefore, the proof is completed.

\end{proof}

		\subsection{Superconvergence}In this subsection, we show the superconvergence for $\|\bm{Q}_h\bm{u}-\bm{u}_h\|_{L^2(\Omega)}$, for which we shall apply the duality argument. To this end, 
we 	consider the dual problem
				\begin{alignat}{2}
					\tdiv \underline{\psi}&=\bm{Q}_h\bm{u}-\bm{u}_h&& \quad \mbox{in}\;\Omega\label{eq:dual1},\\
				\mathcal{A}	\underline{\psi}&=-\varepsilon( \bm{\phi})&&\quad \mbox{in}\;\Omega\label{eq:dual2},
				\end{alignat}
				which satisfies the following elliptic regularity estimate
				\begin{align}
					\|\bm{\phi}\|_{H^2(\Omega)}\leq C \|\bm{Q}_h\bm{u}-\bm{u}_h\|_{L^2(\Omega)}.\label{eq:dual}
				\end{align}
				This estimate is known to hold, for instance, if the domain $\Omega$ is convex (cf. \cite{Grisvard11}).
				
In addition, we let $\mathcal{A}\underline{\psi}=-\nabla \bm{\phi}+\zeta \underline{\delta}$, where $\zeta = \text{rot}(\bm{\phi})$.  
%
We let $\bm{\phi}_h$, $\underline{\psi}_h$ and $\zeta_h$ denote the corresponding approximations of $\bm{\phi}$, $\underline{\psi}$ and $\zeta$ based on the discrete formulation \eqref{eq:MPFA1}-\eqref{eq:MPFA3} with the right-hand side $\bm{f}$ replaced by $\bm{Q}_h\bm{u}-\bm{u}_h$. Then it holds in view of Theorem~\ref{thm:uL2} and \eqref{eq:dual}
					\begin{equation}
						\begin{split}
					\| \underline{\psi}-\underline{\psi}_h\|_{L^2(\Omega)}+\|\zeta-\zeta_h\|_{L^2(\Omega)}&\leq C \Big( h \|\underline{\psi}\|_{H^1(\Omega)}+h\|\zeta\|_{H^1(\Omega)}+h\|\bm{Q}_h\bm{u}-\bm{u}_h\|_{L^2(\Omega)} \Big)\\
					&\leq C h\|\bm{Q}_h\bm{u}-\bm{u}_h\|_{L^2(\Omega)}.
					\end{split}
					\label{eq:con-dual}
				\end{equation}

		In the next theorem, we show the superconvergence, for which we assume that the solution is smooth enough.
				
				\begin{theorem}(superconvergence).
					Let $(\underline{\sigma},\bm{u})$ be the exact solution and let $(\underline{\sigma}_h,\bm{u}_h,\gamma_h)\in \underline{\Sigma}_h\times \bm{U}_h\times \Gamma_h$ be the discrete solution of \eqref{eq:MPFA1}-\eqref{eq:MPFA3} obtained by either Method 1 or Method 2. Assume that $\underline{\sigma}\in \underline{H}^1(\Omega)$ with $\tdiv \underline{\sigma}\in \bm{H}^1(\Omega)$ and $\gamma\in H^1(\Omega)$, then the following convergence error estimate holds
					\begin{align*}
						\|\bm{Q}_h\bm{u}-\bm{u}_h\|_{L^2(\Omega)}\leq C h^2\Big(\|\underline{\sigma}\|_{H^1(\Omega)}+\|\gamma\|_{H^1(\Omega)}+\|\tdiv\underline{\sigma}\|_{H^1(\Omega)}\Big).
					\end{align*}
				\end{theorem}

				\begin{proof}
					We multiply \eqref{eq:dual1} by $\underline{\sigma}-\underline{\sigma}_h$ and \eqref{eq:dual2} by $\bm{Q}_h\bm{u}-\bm{u}_h$, and using the fact that  $\mathcal{A}\underline{\psi}=-\nabla \bm{\phi}+\zeta \underline{\delta}$, we have
					\begin{equation}
						\begin{split}
						\|\bm{Q}_h\bm{u}-\bm{u}_h\|_{L^2(\Omega)}^2&=(\bm{Q}_h\bm{u}-\bm{u}_h, \tdiv \underline{\psi})+(\mathcal{A}\underline{\psi},\underline{\sigma}-\underline{\sigma}_h)+(\nabla \bm{\phi},\underline{ \sigma}-\underline{\sigma}_h)\\
						&-(\zeta, \text{as}(\underline{\sigma}-\underline{\sigma}_h))+(\text{as}(\underline{\psi}),\gamma-\gamma_h).\end{split}\label{eq:Qhpph}
					\end{equation}
			In view of the definition of $\underline{\Pi}_h$, we can obtain
				\begin{align}
				(\tdiv \underline{\Pi}_h(\underline{\sigma}-\underline{\sigma}_h),\bm{v})&=0\quad \forall \bm{v}\in \bm{U}_h\label{eq:error2-RT}.
				\end{align}
				Taking $\bm{v}=\tdiv \underline{\Pi}_h(\underline{\sigma}-\underline{\sigma}_h)$ in \eqref{eq:error2-RT} above, we have
			\begin{align*}
			\tdiv \underline{\Pi}_h(\underline{\sigma}-\underline{\sigma}_h)=0,
			\end{align*}
			which yields
			\begin{align}
				\tdiv (\underline{\sigma}-\underline{\Pi}_h\underline{\sigma}_h)=\tdiv (\underline{\sigma}-\underline{\Pi}_h\underline{\sigma}).\label{eq:divuuh}
			\end{align}
					It follows from \eqref{eq:error1-RT1}-\eqref{eq:error3-RT} by setting $\underline{w}=\underline{\psi}_h$, $\bm{v}=\bm{\phi}_h$ and $\xi = \zeta_h$ that
					\begin{align}
						(\mathcal{A}(\underline{\sigma}-\underline{\sigma}_h),\underline{\psi}_h)+(\bm{u}-\bm{u}_h,\tdiv\underline{\Pi}_h \underline{\psi}_h)+(\gamma-\gamma_h, \text{as}(\underline{\psi}_h))&=(\mathcal{A}\underline{\sigma}+\gamma\underline{\delta},\underline{\psi}_h-\underline{\Pi}_h\underline{\psi}_h),\label{eq:psi1}\\
						(\tdiv\underline{\Pi}_h (\underline{\sigma}-\underline{\sigma}_h),\bm{\phi}_h)&=0.\label{eq:psi2}\\
						(\textnormal{as}(\underline{\sigma}-\underline{\sigma}_h), \zeta_h)&=0.\label{eq:psi3}
					\end{align}
					As a consequence of \eqref{eq:psi1}-\eqref{eq:psi3} and \eqref{eq:Qhpph}, we have
					\begin{equation}
						\begin{split}
							\|\bm{Q}_h\bm{u}-\bm{u}_h\|_{L^2(\Omega)}^2&=(\bm{Q}_h\bm{u}-\bm{u}_h, \tdiv \underline{\psi})+(\mathcal{A}\underline{\psi},\underline{\sigma}-\underline{\sigma}_h)+(\nabla \bm{\phi},\underline{ \sigma}-\underline{\sigma}_h)-(\zeta-\zeta_h, \text{as}(\underline{\sigma}-\underline{\sigma}_h))\\
							&+(\text{as}(\underline{\psi}),\gamma-\gamma_h)\\
							&=(\bm{Q}_h\bm{u}-\bm{u}_h, \tdiv (\underline{\psi}-\underline{\Pi}_h\underline{\psi}_h))+(\mathcal{A}(\underline{\sigma}-\underline{\sigma}_h),\underline{\psi}-\underline{\psi}_h)+(\nabla \bm{\phi},\underline{ \sigma}-\underline{\sigma}_h)\\
							&+(\mathcal{A}\underline{\sigma}+\gamma\underline{\delta},\underline{\psi}_h-\underline{\Pi}_h\underline{\psi}_h)-(\zeta-\zeta_h, \text{as}(\underline{\sigma}-\underline{\sigma}_h))+(\gamma-\gamma_h, \text{as}(\underline{\psi}-\underline{\psi}_h)).
						\end{split}\label{eq:Qhpdual}
					\end{equation}
					As $\underline{\psi}_h$ is the approximation of $\underline{\psi}$, it holds by proceeding similarly to \eqref{eq:error2-RT}
					\begin{align*}
						(\bm{Q}_h\bm{u}-\bm{u}_h,\tdiv (\underline{\psi}-\underline{\Pi}_h\underline{\psi}_h))=0.
					\end{align*}
					The Cauchy-Schwarz inequality and \eqref{thm:uL2} lead to
					\begin{align*}
						(\mathcal{A}(\underline{\sigma}-\underline{\sigma}_h),\underline{\psi}-\underline{\psi}_h)&\leq \|\underline{\sigma}-\underline{\sigma}_h\|_{\mathcal{A}}\|\underline{\psi}-\underline{\psi}_h\|_{L^2(\Omega)}\leq Ch^2\Big( \|\underline{\sigma}\|_{H^1(\Omega)}+\|\gamma\|_{H^1(\Omega)}\Big) \|\bm{Q}_h\bm{u}-\bm{u}_h\|_{L^2(\Omega)}.
					\end{align*}
					Now we estimate the third term on the right-hand side of \eqref{eq:Qhpdual}. To this end, an application of integration by parts implies
					\begin{align*}
						(\nabla \bm{\phi}, \underline{\sigma}-\underline{\sigma}_h)&=(\nabla \bm{\phi}, \underline{\sigma}-\underline{\Pi}_h \underline{\sigma}_h)+(\nabla \bm{\phi},\underline{\Pi}_h\underline{\sigma}_h-\underline{\sigma}_h)\\
						&=-(\bm{\phi}, \tdiv (\underline{\sigma}-\underline{\Pi}_h\underline{\sigma}_h))-(\underline{\Pi}_h\underline{\sigma}_h-\underline{\sigma}_h,\mathcal{A}\underline{\psi}-\zeta\underline{\delta}),
					\end{align*}
					which coupled with \eqref{eq:psi2} and \eqref{eq:divuuh} yields
					\begin{align*}
						(\nabla \bm{\phi}, \underline{\sigma}-\underline{\sigma}_h)
						&=-(\bm{\phi}-\bm{\phi}_h, \tdiv(\underline{\sigma}-\underline{\Pi}_h\underline{\sigma}_h))-(\underline{\Pi}_h\underline{\sigma}_h-\underline{\sigma}_h,\mathcal{A}\underline{\psi}-\zeta\underline{\delta})\\
						&=-(\bm{\phi}-\bm{\phi}_h, \tdiv(\underline{\sigma}-\underline{\Pi}_h\underline{\sigma}))-(\underline{\Pi}_h\underline{\sigma}_h-\underline{\sigma}_h,\mathcal{A}\underline{\psi}-\zeta\underline{\delta}).
					\end{align*}
					The first term on the right-hand side is bounded by the Cauchy-Schwarz inequality and the interpolation error estimate
					\begin{align*}
						(\bm{\phi}-\bm{\phi}_h, \tdiv(\underline{\sigma}-\underline{\Pi}_h\underline{\sigma}))\leq C h^2 \|\bm{\phi}\|_{H^2}\|\tdiv\underline{\sigma}\|_{H^1(\Omega)}.
					\end{align*}
					It follows from Lemma~\ref{lemma:phi} that
					\begin{align*}
						(\underline{\Pi}_h\underline{\sigma}_h-\underline{\sigma}_h,\mathcal{A}\underline{\psi}-\zeta\underline{\delta})&\leq C h(\|\underline{\psi}\|_{H^1(\Omega)}+\|\zeta\|_{H^1(\Omega)}) \|\underline{\Pi}_h\underline{\sigma}_h-\underline{\sigma}_h\|_{L^2(\Omega)}.
					\end{align*}
					The triangle inequality yields
					\begin{align}
						\|\underline{\Pi}_h\underline{\sigma}_h-\underline{\sigma}_h\|_{L^2(\Omega)}\leq \|\underline{\Pi}_h(\underline{\sigma}_h-\underline{\sigma})\|_{L^2(\Omega)}+\|\underline{\Pi}_h\underline{\sigma}-\underline{\sigma}\|_{L^2(\Omega)}+ \|\underline{\sigma}-\underline{\sigma}_h\|_{L^2(\Omega)}.\label{eq:Pihuuh}
					\end{align}
					It holds in view of \eqref{eq:Pihdef}
					\begin{align*}
						(\underline{\Pi}_h(\underline{\sigma}-\underline{\sigma}_h)\bm{n},\bm{\mu})_e=((\underline{\Pi}_h\underline{\sigma}-\underline{\sigma}_h)\bm{n}, \bm{\mu})_e\quad \forall\bm{\mu}\in \bm{P}_0(e), e\in \mathcal{F}_{pr}.
					\end{align*}
					Then the equivalence of norms on finite dimensional spaces yields
					\begin{align*}
						\|	\underline{\Pi}_h(\underline{\sigma}-\underline{\sigma}_h)\|_{L^2(\Omega)}\leq C \|\underline{\Pi}_h\underline{\sigma}-\underline{\sigma}_h\|_{L^2(\Omega)}.
					\end{align*}
					This, \eqref{eq:Pihuuh}, the interpolation error estimate and Theorem~\ref{thm:uL2} yield
					\begin{align*}
						\|\underline{\Pi}_h\underline{\sigma}_h-\underline{\sigma}_h\|_{L^2(\Omega)}\leq C( h\|\underline{\sigma}\|_{H^1(\Omega)}+h\|\gamma\|_{H^1(\Omega)}).
					\end{align*}
					Similarly, the fourth term on the right-hand side of \eqref{eq:Qhpdual} can be estimated by Lemma~\ref{lemma:phi}
					\begin{align*}
						(\mathcal{A}\underline{\sigma}+\gamma\underline{\delta},\underline{\psi}_h-\underline{\Pi}_h\underline{\psi}_h)&\leq C h(\| \underline{\sigma}\|_{H^1(\Omega)}+\|\gamma\|_{H^1(\Omega)})\|\underline{\psi}_h-\underline{\Pi}_h\underline{\psi}_h\|_{L^2(\Omega)}\\
						&\leq Ch^2( \|\underline{\sigma}\|_{H^1(\Omega)}+\|\gamma\|_{H^1(\Omega)} )\|\bm{Q}_h\bm{u}-\bm{u}_h\|_{L^2(\Omega)}.
					\end{align*}
				The last two terms on the right-hand side of \eqref{eq:Qhpdual} can be bounded by the Cauchy-Schwarz inequality, Theorem~\ref{thm:uL2} and \eqref{eq:con-dual}. The proof is completed by using the above estimates.

				\end{proof}
				
				\begin{remark}
					Since $\bm{Q}_h\bm{u}$ is $\mathcal{O}(h^2)$-close to $\bm{u}$ at the center of mass of each element, the above theorem implies that
					\begin{align*}
						\|\bm{u}-\bm{u}_h\|_{0,h}\leq C h^2,
					\end{align*}
					where $\|\bm{u}-\bm{u}_h\|_{0,h}=\Big(\sum_{M\in \mathcal{T}_M}|M| (\bm{u}(m_M)-\bm{u}_h)^2\Big)^{1/2}$ and $m_M$ is the center of mass of $M\in \mathcal{T}_M$.
				\end{remark}

\subsection{Reduction to a cell-centered system}\label{sec:cell-center}

The algebraic system yielded by \eqref{eq:MPFA1}-\eqref{eq:MPFA3} has the following form 
\begin{equation}
	\bm{u}=\begin{pmatrix}
	A_{\underline{\sigma} \underline{\sigma}} & A^T_{\underline{\sigma}\bm{u}} & A^T_{\underline{\sigma}\gamma}\\
	A_{\underline{\sigma}\bm{u}}& 0& 0 \\
	A_{\underline{\sigma}\gamma}&0 &0
	\end{pmatrix} \begin{pmatrix}
	\underline{\sigma} \\
	\bm{u} \\
	\gamma
\end{pmatrix} =
\begin{pmatrix}
	0 \\
	F \\
0
\end{pmatrix} \label{eq:m1}
\end{equation}
where 	$(A_{\underline{\sigma} \underline{\sigma}})_{ij}=(\mathcal{A}\underline{w}_j,\underline{w}_i)$,
	$(A_{\underline{\sigma}\bm{u}})_{ij}=\sum_{e\in \mathcal{F}_{dl}}(\bm{v}_j,\jump{\underline{w}_i\bm{n}})_e$, and
	$(A_{\underline{\sigma}\gamma})_{ij}=(\text{as}(\underline{\sigma}_j), \xi_i)$. The definitions of 
	stress basis functions $\underline{w}_i$  result in matrix $A_{\underline{\sigma} \underline{\sigma}}$ being block-diagonal with 
	blocks associated with the mesh vertices. Specifically, for a vertex $p$ shared by $k$ edges for faces $e_1,\cdots,e_k$ ($k=4$ in 2\text{D} and $k=12$ in 3\text{D}) as shown in Figure \ref{fig:dof-MPFA}. Let $\underline{w}_1,\underline{w}_2\cdots \underline{w}_{d},\cdots \underline{w}_{dk}$ be the stress basis functions associated with the vertex $p$, based on the definition of $\underline{w}_{i}$, $(\mathcal{A}\cdot,\cdot)$ localizes  basis functions interaction
	around mesh vertices by decoupling them from the rest of the basis functions, as a result, the matrix $A_{\underline{\sigma} \underline{\sigma}}$ is block-diagonal with $dk\times dk$ blocks associated with mesh vertices. It is obvious the bilinear form $(\mathcal{A}\underline{\sigma},\underline{w})$ is 
	an inner product on $\underline{\Sigma}_h$  and $(\mathcal{A}{\underline{\sigma}},\underline{\sigma})^{1/2}$ is a norm in $\underline{\Sigma}_h$ equivalent 
	to $\|\cdot \|_{L^2(\Omega)}$, therefore, the blocks of $A_{\underline{\sigma} \underline{\sigma}}$ are symmetric and positive definite. 
	 Therefore, the  stress $\underline{\sigma}$ can be
	easily eliminated from \eqref{eq:m1} by inverting small local linear systems, which yields a displacement-rotation system 
	
\begin{equation}
\begin{pmatrix}
A_{\underline{\sigma}\bm{u}}A_{\underline{\sigma} \underline{\sigma}}^{-1}A_{\underline{\sigma}\bm{u}}^T  & A_{\underline{\sigma}\bm{u}}A_{\underline{\sigma} \underline{\sigma}}^{-1}A_{\underline{\sigma}\gamma}^T\\
A_{\underline{\sigma}\gamma}A_{\underline{\sigma} \underline{\sigma}}^{-1}A_{\underline{\sigma}\bm{u}}^T  & A_{\underline{\sigma}\gamma}A_{\underline{\sigma} \underline{\sigma}}^{-1}A_{\underline{\sigma}\gamma}^T \\
\end{pmatrix} 
\begin{pmatrix}
	\bm{u} \\
	\gamma \\
\end{pmatrix} 	
=
\begin{pmatrix}
	\widetilde{F} \\
	\widetilde{H} \\
\end{pmatrix} 	\label{eq:m2}
\end{equation}
Clearly, the matrix in \eqref{eq:m2} is symmetric, the inf-sup condition \eqref{eq:inf} implies it is also positive definite. Indeed, we can state the following lemma.
\begin{lemma} The displacement-rotation system \eqref{eq:m2} is symmetric and positive definite.

\end{lemma}

\begin{proof}
	We only need to prove the positive definiteness. To this end, for any $(\bm{v}^T\quad  \bm{\xi}^T)\neq 0$, we have
	\begin{align*}
		(\bm{v}^T\quad  \bm{\xi}^T)\begin{pmatrix}
			A_{\underline{\sigma}\bm{u}}A_{\underline{\sigma} \underline{\sigma}}^{-1}A_{\underline{\sigma}\bm{u}}^T  & A_{\underline{\sigma}\bm{u}}A_{\underline{\sigma} \underline{\sigma}}^{-1}A_{\underline{\sigma}\gamma}^T\\
			A_{\underline{\sigma}\gamma}A_{\underline{\sigma} \underline{\sigma}}^{-1}A_{\underline{\sigma}\bm{u}}^T  & A_{\underline{\sigma}\gamma}A_{\underline{\sigma} \underline{\sigma}}^{-1}A_{\underline{\sigma}\gamma}^T \\
		\end{pmatrix} 
	\begin{pmatrix}
		\bm{v}\\\bm{\xi}
	\end{pmatrix}=(A_{\underline{\sigma}\bm{u}}^T\bm{v}+A_{\underline{\sigma}\gamma}^T\bm{\xi})^T A_{\underline{\sigma} \underline{\sigma}}^{-1}(A_{\underline{\sigma}\bm{u}}^T\bm{v}+A_{\underline{\sigma}\gamma}^T\bm{\xi})>0
	\end{align*}
owing to the inf-sup condition \eqref{eq:inf}.

\end{proof}
For Method 2, Equation \eqref{eq:m2} can not be reduced further since the displacement and rotation are coupled 
and thus it is the final linear system we need to solve.  However, for Method 1, we can 
further eliminate the rotation. More precisely, the stress basis and rotation basis associated with a vertex have no interaction with stress basis and rotation basis associated with an another vertex, which  implies that 
$A_{\underline{\sigma}\gamma}$ is block-diagonal with $d(d-1)/2\times dk$ blocks, which means the rotation matrix 
$A_{\underline{\sigma}\gamma}A_{\underline{\sigma} \underline{\sigma}}^{-1}A_{\underline{\sigma}\gamma}^T$ is also block-diagonal with $d(d-1)/2\times d(d-1)/2$ blocks and thus can be easily inverted. 
After elimination of rotation, Equation \eqref{eq:m2} is reduced to 
\begin{equation}
	\big(A_{\underline{\sigma}\bm{u}}A_{\underline{\sigma} \underline{\sigma}}^{-1}A_{\underline{\sigma}\bm{u}}^T-
	A_{\underline{\sigma}\bm{u}}A_{\underline{\sigma} \underline{\sigma}}^{-1}A_{\underline{\sigma}\gamma}^T
	(A_{\underline{\sigma}\gamma}A_{\underline{\sigma} \underline{\sigma}}^{-1}A_{\underline{\sigma}\gamma}^T)^{-1}
	A_{\underline{\sigma}\gamma}A_{\underline{\sigma} \underline{\sigma}}^{-1}A_{\underline{\sigma}\bm{u}}^T\big)\bm{u}=\widehat{F}.
	\label{eq:m3}
\end{equation}
 The matrix in \eqref{eq:m3} is symmetric and positive definite, since it is a Schur complement of the symmetric
 and positive definite matrix in \eqref{eq:m2}; see \cite{HornJohnson2012}. We remark that one can easily design multigrid preconditioners to efficiently solve \eqref{eq:m2} and 
 \eqref{eq:m3}.

\section{Extension to Darcy flow and the Stokes equations}\label{sec:extension}
In this section, we show the extension of the proposed scheme for Darcy flow and the Stokes equations. 
\subsection{Darcy flow}
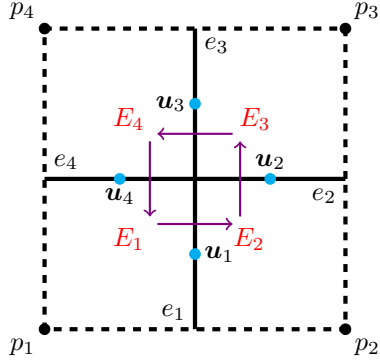
\begin{figure}[t]
	\begin{center}
		\begin{tikzpicture}[scale=2]
			\coordinate (O) at (0,0);
			\coordinate (A) at (1,0);
			\coordinate (B) at (1,1);
			\coordinate (C) at (0,1);
			\coordinate (D) at (-1,1);
			\coordinate (E) at (-1,0);
			\coordinate (F) at (-1,-1);
			\coordinate (G) at (0,-1);
			\coordinate (H) at (1,-1);
			
			\draw[fill, black] (B) circle (1pt);
			\draw[fill, black] (D) circle (1pt);
			\draw[fill, black] (F) circle (1pt);
			\draw[fill, black] (H) circle (1pt);
			
			\draw[ultra thick, dashed] (A) -- (B) -- (C) -- (D) -- (E) -- (F) -- (G) -- (H) -- cycle;
			\draw[ultra thick] (A) -- (O) -- (E);
			\draw[ultra thick] (C) -- (O) -- (G);
			
			\coordinate (LLC1) at ($0.5*(O)+0.5*(A)$);
			\coordinate (LLC2) at ($0.5*(O)+0.5*(C)$);
			\coordinate (LLC3) at ($0.5*(O)+0.5*(E)$);
			\coordinate (LLC4) at ($0.5*(O)+0.5*(G)$);
			
			\draw[violet,->,thick] ($(LLC1)+(-0.2,-0.25)$) -- ($(LLC1)+(-0.2,0.25)$);
			\draw[violet,->,thick] ($(LLC2)+(0.25,-0.2)$) -- ($(LLC2)+(-0.25,-0.2)$);
			\draw[violet,->,thick] ($(LLC3)+(0.2,0.25)$) -- ($(LLC3)+(0.2,-0.25)$);
			\draw[violet,->,thick] ($(LLC4)+(-0.25,0.2)$) -- ($(LLC4)+(0.25,0.2)$);
			
			\draw[fill, cyan] (LLC1) circle (1pt);
			\draw[fill, cyan] (LLC2) circle (1pt);
			\draw[fill, cyan] (LLC3) circle (1pt);
			\draw[fill, cyan] (LLC4) circle (1pt);

			\node[red] at (AC1) {$E_3$};
			\node[red] at (AC2) {$E_4$};
			\node[red] at (AC3) {$E_1$};
			\node[red] at (AC4) {$E_2$};
			
			\node[above right] at (B) {$p_3$};
			\node[above left] at (D) {$p_4$};
			\node[below left] at (F) {$p_1$};
			\node[below right] at (H) {$p_2$};
			
			\node[above] at (LLC1) {$\bm{u}_2$};
			\node[left] at (LLC2) {$\bm{u}_3$};
			\node[below] at (LLC3) {$\bm{u}_4$};
			\node[right] at (LLC4) {$\bm{u}_1$};
			
			\node[below left] at (A) {$e_2$};
			\node[below right] at (C) {$e_3$};
			\node[above right] at (E) {$e_4$};
			\node[above left] at (G) {$e_1$};
			
		\end{tikzpicture}
	\end{center}
	\caption{One interaction region with four subcells numbered $E_i, i=1,\cdots,4$, in the reference space. The filled dots denote the cell pressure $\{p_i\}$ and the blue dots denote the velocity $\{\bm{u}_i\}$.}
	\label{fig:dof-intersectionp}
\end{figure}

We consider the following Darcy flow
\begin{alignat*}{2}
	\bm{u}&=-\nabla p&&\quad \mbox{in}\;\Omega,\\
	\nabla\cdot \bm{u}&=f&&\quad \mbox{in}\;\Omega.
\end{alignat*}
The discrete formulation reads as follows: Find $(\bm{u}_h,p_h)\in \bm{U}_h\times \Gamma_h^2$ such that
\begin{align}
	(\bm{u}_h,\bm{v})+\sum_{e\in \mathcal{F}_{dl}} (p_h,\jump{\bm{v}\cdot\bm{n}})_e&=0\quad \forall \bm{v}\in \bm{U}_h,\label{eq:MPFA1D}\\
	\sum_{e\in \mathcal{F}_{pr}^{\frac{1}{2}}}(\bm{u}_h\cdot\bm{n},\jump{q})_e&=(f,q)\quad \forall q\in \Gamma_h^2.\label{eq:MPFA2D}
\end{align}
The mass matrix corresponding to the first term of \eqref{eq:MPFA1D} is block diagonal, and we can eliminate $\bm{u}_h$ using \eqref{eq:MPFA1D}. More specifically, $\bm{u}_h$ can be expressed in terms of $p_h$ over each interaction region, that is,
\begin{align*}
	\bm{u}_h=A^{-1} b,
\end{align*} 
where $A$ is the mass matrix for $\bm{u}_h$ corresponding to each interaction region and $p=-(p_2-p_1,p_3-p_2,p_3-p_4,p_4-p_1)$; see Figure~\ref{fig:dof-intersectionp}. We can observe that this elimination is in a similar fashion to that of \cite[3.13]{Klausen06}.

\subsection{The Stokes equations}

In this subsection, we consider the Stokes equations given by
\begin{align}
	-\Delta \bm{u}+\nabla p =\bm{f} \quad \mbox{in}\;\Omega,\label{eq:Stokes1}\\
	\nabla\cdot \bm{u}=0\quad \mbox{in}\;\Omega,\label{eq:Stokes2}
\end{align}
which are supplemented with homogeneous Dirichlet boundary condition. 

We introduce the auxiliary unknown $\underline{w}=\nabla \bm{u}$, then \eqref{eq:Stokes1}-\eqref{eq:Stokes2} can be recast into the following equivalent form
\begin{align*}
	\underline{w}&=\nabla \bm{u} \quad \mbox{in}\;\Omega,\\
	-\tdiv \underline{w}+\nabla p&=\bm{f}\quad \mbox{in}\;\Omega,\\
	\nabla\cdot \bm{u}&=0\quad \mbox{in}\;\Omega.
\end{align*}
We define the following bilinear forms
\begin{align*}
	B_h(\underline{w},\bm{v})&=-\sum_{e\in \mathcal{F}_{pr}^{\frac{1}{2}}}(\underline{w}\bm{n},\jump{\bm{v}})_e,
	\\B_h^*(\bm{v},\underline{w})&=\sum_{e\in \mathcal{F}_{dl}}(\bm{v},\jump{\underline{w}\bm{n}})_e
\end{align*}
and
\begin{align*}
	b_h^*(q,\bm{v})&=\sum_{e\in \mathcal{F}_{pr}^{\frac{1}{2}}}(q,\jump{\bm{v}\cdot\bm{n}})_e,\\
	b_h(\bm{v},q)&=-\sum_{e\in \mathcal{F}_{dl}} (\bm{v}\cdot\bm{n},\jump{q})_e.
\end{align*}
Then the discrete formulation for the Stokes equations reads as follows: Find $(\underline{w}_h,\bm{u}_h,p_h)\in \underline{\Sigma}_h\times \bm{U}_h\times \Gamma_h$ such that
\begin{alignat}{2}
	B_h(\underline{w}_h,\bm{v})+b_h^*(p_h,\bm{v})&=(\bm{f},\bm{v}) &&\quad \forall \bm{v}\in \bm{U}_h\label{eq:discrete1}\\
	(\underline{w}_h,\underline{H})&=B_h^*(\bm{u}_h,\underline{H}) &&\quad \forall \underline{H}\in \underline{\Sigma}_h,\label{eq:discrete2}\\
	b_h(\bm{u}_h,q)&=0&&\quad \forall q\in \Gamma_h.\label{eq:discrete3}
\end{alignat}
In addition, the following adjoint properties hold
\begin{align*}
	B_h(\underline{w},\bm{v})&=B_h^*(\bm{v},\underline{w}) \quad \forall (\underline{w},\bm{v})\in \underline{\Sigma}_h\times \bm{U}_h,\\
	b_h(\bm{v},q)&=b_h^*(q,\bm{v})\quad \forall (\bm{v},q)
	\in \bm{U}_h\times \Gamma_h.
\end{align*}
Then one can perform the error analysis in a similar fashion to that of the linear elasticity problem. We can also apply local elimination for the velocity gradient to improve the computational efficiency.

\section{Numerical experiments}\label{sec:numerical}
		In this section, we present several examples to verify the convergence of our method.  Specifically, the relative errors of displacement $\bm{u}_h$, stress tensor $\underline{\sigma}_h$ and rotation $\gamma_h$ and the corresponding convergence rates will be reported.
	For simplicity of notation, we denote $\|a-b\|_{L^2(\Omega)}$ as the relative $L^2$-error  between $a$ and $b$ instead of using $\|a-b\|_{L^2(\Omega)}/\|b\|_{L^2(\Omega)}$. Moreover, $\underline{\sigma}_h$ has different values in subcells of an macro-element, we compute the mean value of  $\underline{\sigma}_h$ over  a
	macro-element denoted by $\mathcal{M}\underline{\sigma}_h$ and its relative error against exact solution in the center of each macro-element.  

	\textbf{Data availability}: The source code for generating the data in all tables can be found in https://github.com/aggietx/a-new-mixed-fem.
	
\begin{figure}[htb!]
		\centering
		\includegraphics[trim={0cm .8cm 0cm 0cm},clip,width=2in]{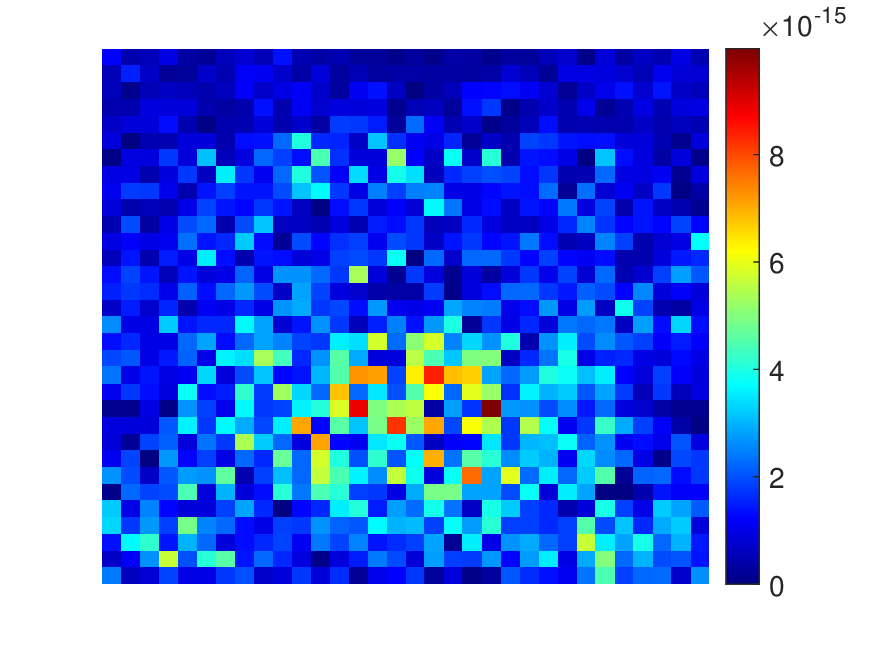}			
		\caption{\label{fig:conser}residual distribution, $h=1/32$, example 1.}
		
	\end{figure}

	\hong{\subsection{Example 1: 2D structured grid case}}
	In the first example, the analytical solution is given by 
	
	\begin{equation*}
		\bm{u}=\begin{pmatrix}
			\cos(\pi x)\sin(2\pi y)\\
			\cos(\pi y)\sin(\pi x)
		\end{pmatrix} 
	\end{equation*}
	with Dirichlet boundary condition, and the body force is then computed using Lam\'e coefficients $\lambda = 123, \mu = 79.3$.
	Table \ref{ex1} shows errors and convergence rates  on a sequence of uniform mesh refinements, computed using rotation over the interaction region (Method 1) and macro-element (Method 2). 
	Convergence  rates of $\bm{u}_h$ and $\underline{\sigma}_h$ agree with the theoretical error analysis. Superconvergence of $\mathcal{M}\underline{\sigma}_h$ and rotation $\gamma_h$ are observed. 
	Figure    \ref{fig:conser} shows an example of residual distribution which clearly demonstrates the local conservation property (see Remark \ref{re:local}) of the proposed method.
	
	\hong{\subsection{Example 2: 3D structured grid case}}
	
	The second example aims to illustrate the convergence behavior of our methods in three dimension, in particular, the analytical solution is 
	
	\begin{equation*}
		\bm{u}=\begin{pmatrix}
			0\\
			-(e^x-1)(y-\cos(\frac{\pi}{12})(y-\frac{1}{2})+\sin(\frac{\pi}{12})(z-\frac{1}{2})-\frac{1}{2})\\
			-(e^x-1)(z-\sin(\frac{\pi}{12})(y-\frac{1}{2})-\cos(\frac{\pi}{12})(z-\frac{1}{2})-\frac{1}{2})
		\end{pmatrix} ,
	\end{equation*}
	and Dirichlet boundary condition is considered with $\lambda =\mu = 79.3$. Relative errors and convergence rates for both methods are presented in Table \ref{ex2}.
	It is observed for rotation defined over the macro-element,  convergence rates of all variables are similar as the first example, for rotation defined over the interaction region, 
	convergence  rate of $\gamma_h$ and $\underline{\sigma}_h$ are in accordance with error analysis, there are also  superconvergence of $\mathcal{M}\underline{\sigma}_h$ and rotation $\gamma_h$.

\begin{figure}
\centering
\subfigure[Computed displacement of example 3,  $h=1/48$ ]{
	\includegraphics[trim={4cm 0.3cm 1cm 0.2cm},clip,width=2in]{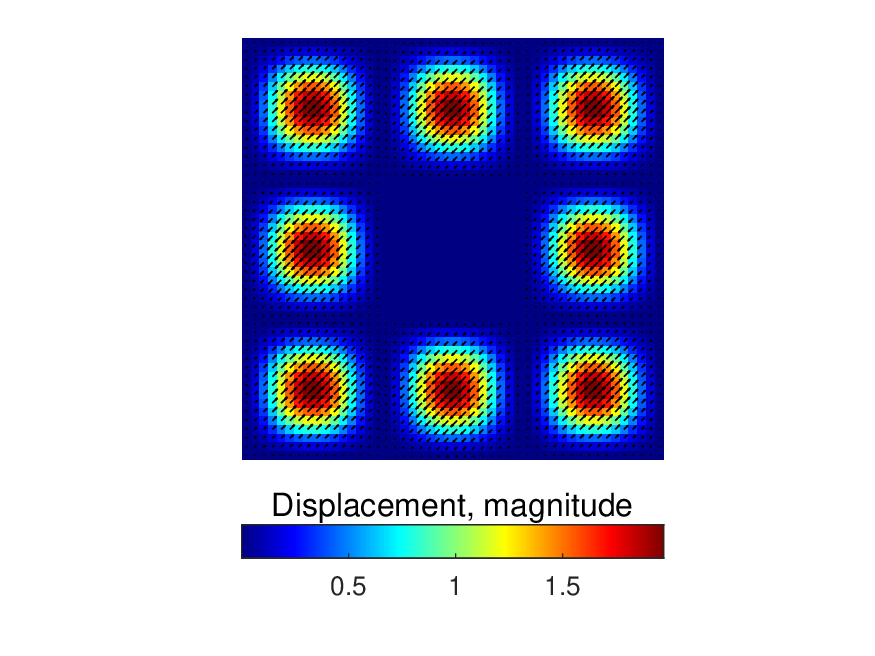}}	
\subfigure[Computed displacement of example 1 on a $h^2$-parallelogram mesh, $h=1/32$]{
	\includegraphics[trim={4cm 0.3cm 1cm 0.2cm},clip,width=2in]{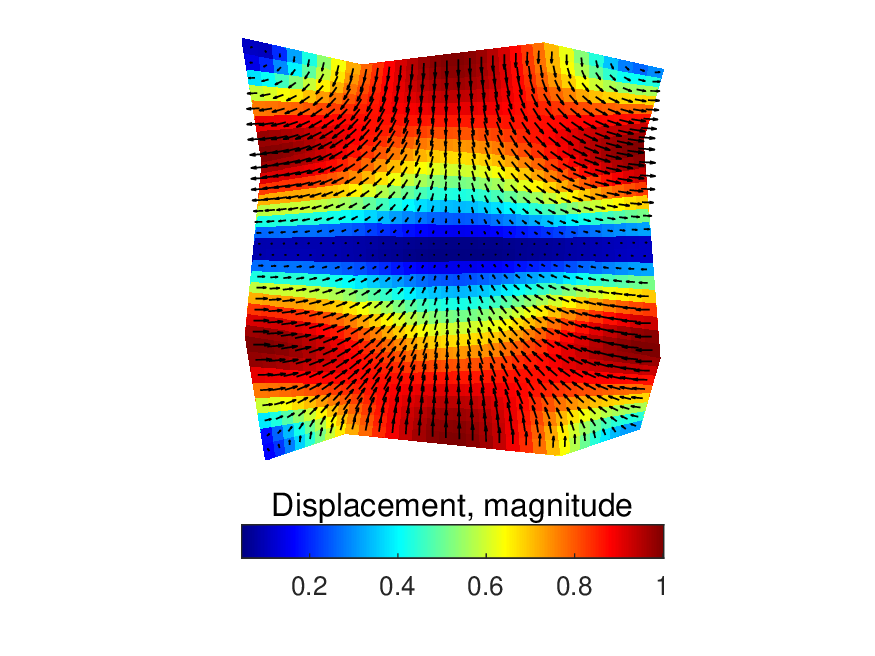}}
\subfigure[Computed displacement  of example 1 on a quadrilateral mesh, $h=1/32$]{
	\includegraphics[trim={4cm 0.3cm 1cm 0.2cm},clip,width=2in]{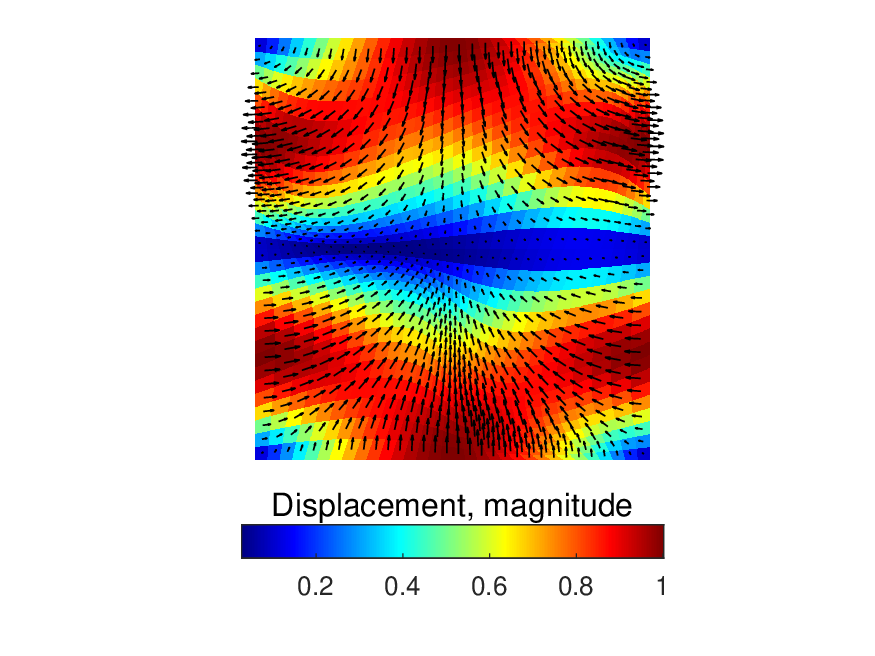}}			
\caption{Computed displacements. }
\label{fig:sol}
\end{figure}

	\hong{\subsection{Example 3: heterogeneous materials }}
	The third example is taken from \cite{Ambartsumyan20simplicial} and we aim to demonstrating the performance of the proposed methods for heterogeneous materials. We set the analytical solution $\bm{u}$ 
	as 
	\begin{equation*}
		\bm{u}=\frac{1}{(1-\chi)+\kappa \chi}\begin{pmatrix}
			\sin(3\pi x)\sin(3\pi y)\\
			\sin(3\pi x)\sin(3\pi y)
		\end{pmatrix} ,
	\end{equation*}
	where $\chi(x,y)$ is defined as 
	\begin{equation*}
		\chi(x,y)=\left\{
		\begin{aligned}
			1 & \quad \text{  if } \text{min}(x,y)>\frac{1}{3} \text{ and } \text{ max}(x,y) <\frac{2}{3}, \\
			0 & \quad \text{  otherwise.}
		\end{aligned}
		\right.
	\end{equation*}
	We set $\kappa=10^6$ and $\lambda =\mu =(1-\chi)+\kappa \chi$. Since the coefficient here is highly heterogeneous, we consider the modified weak formulation \eqref{eq:MPFA1b}-\eqref{eq:MPFA3b} for Method 1.
	Table \ref{ex3} lists errors and convergences rates of this heterogeneous example. Similar convergence rates as example 2 can be found, 
	it is worthwhile to mention superconvergence of $\mathcal{M}\underline{\sigma}_h$ and rotation $\gamma_h$ are also observed. The computed displacement for this example are displayed in the left panel of 
	Figure \ref{fig:sol}.
	\hong{\subsection{Example 4: nearly incompressible materials}}
	Example 4 intends to check the locking-free property of the proposed method, specifically, the exact solution is given by 
	\begin{equation*}
		\bm{u}=\begin{pmatrix}
			\sin(\pi x)\sin(\pi y)+\frac{1}{2\lambda}x\\
			\cos(\pi x)\cos(\pi y)+\frac{1}{2\lambda}y
		\end{pmatrix} ,
	\end{equation*}
	and the load $\bm{f}$ is 
	\begin{equation*}
		\bm{f}=\begin{pmatrix}
			2\pi^2\sin(\pi x)\sin(\pi y)\\
			2\pi^2\cos(\pi x)\cos(\pi y)
		\end{pmatrix} .
	\end{equation*}
	Here, we choose $\mu=1$ and $\lambda=10^6$.  The analytical solution  has vanishing
	divergence in the limit $\lambda \rightarrow +\infty$, besides, $\bm{f}$ does not
	depend on $\lambda$ , therefore, this is an ideal example for checking numerically that our proposed
	method is indeed locking-free. We report errors and convergence rates in Table \ref{ex4}, as expected, no deteriorate of convergence rates are observed which
	clearly demonstrates the locking-free property of the proposed method. 
	
	\begin{table}
		\centering
		\begin{adjustbox}{max width=\textwidth}
			
			\begin{tabular}{ccccccccc}
				\toprule 
				
				\multicolumn{9}{c}{Method 1-rotation on the interaction region}\tabularnewline
				\midrule 
				\multirow{2}{*}{$h$}  & \multicolumn{2}{c}{$\|\underline{\sigma}-\underline{\sigma}_h\|_{L^2(\Omega)}$} & \multicolumn{2}{c}{$\|\underline{\sigma}-\mathcal{M}\underline{\sigma}_h\|_{L^2(\Omega)}$}& \multicolumn{2}{c}{$\|\bm{u}-\bm{u}_h\|_{L^2(\Omega)}$}& \multicolumn{2}{c}{$\|\gamma-\gamma_h\|_{L^2(\Omega)}$}\tabularnewline
				\cmidrule{2-9} \cmidrule{3-9} \cmidrule{4-9} \cmidrule{5-9} \cmidrule{6-9} \cmidrule{7-9} \cmidrule{8-9} \cmidrule{9-9} 
				& Error & Rate & Error & Rate &Error & Rate&Error & Rate\tabularnewline
				\midrule 
				$1/4$& 1.9615E-01 &   /&8.4922E-02& /&  1.1917E-01 &  /&1.4999E-01 & /\tabularnewline\midrule
				$1/8$  &1.0045E-01  &   0.9655&2.6872E-02& 1.6600& 2.8380E-02  &2.0701& 4.4583E-02&1.7503  \tabularnewline\midrule
				$1/16$	 &4.8951E-02  &  1.0371 &6.9991E-03 & 1.9409& 6.9959E-03  &2.0203&1.1770E-02&   1.9214\tabularnewline\midrule
				$1/32$ & 2.4525E-02 &  0.9971 &1.7713E-03 &1.9824&  1.7429E-03 & 2.0050&2.9897E-03&1.9770    \tabularnewline\midrule
				$1/64$     &   1.2262E-02 &1.0001&4.4436E-04 &1.9950& 4.3534E-04&2.0013&7.5065E-04& 1.9938 \tabularnewline\midrule
				$1/128$	 & 6.1292E-03 &  1.0004 &1.1120E-04 &1.9986& 1.0881E-04  &2.0003&1.8787E-04&  1.9984 \tabularnewline
				\toprule 	
				\multicolumn{9}{c}{Method 2-rotation on the macro-element}\tabularnewline	
				\midrule 
				\multirow{2}{*}{$h$}  & \multicolumn{2}{c}{$\|\underline{\sigma}-\underline{\sigma}_h\|_{L^2(\Omega)}$} & \multicolumn{2}{c}{$\|\underline{\sigma}-\mathcal{M}\underline{\sigma}_h\|_{L^2(\Omega)}$}& \multicolumn{2}{c}{$\|\bm{u}-\bm{u}_h\|_{L^2(\Omega)}$}& \multicolumn{2}{c}{$\|\gamma-\gamma_h\|_{L^2(\Omega)}$}\tabularnewline
				\cmidrule{2-9} \cmidrule{3-9} \cmidrule{4-9} \cmidrule{5-9} \cmidrule{6-9} \cmidrule{7-9} \cmidrule{8-9} \cmidrule{9-9}
				& Error & Rate & Error & Rate &Error & Rate&Error & Rate\tabularnewline
				\midrule 
				$1/4$& 1.5018E-01 &  / &6.9032E-02 & /&  1.0630E-01 & /&9.8917E-02&   /\tabularnewline	\midrule 
				$1/8$&6.9372E-02  & 1.1143  &1.7791E-02 &1.9561& 2.5806E-02  &  2.0424&1.4475E-02& 2.7727  \tabularnewline	\midrule 
				$1/16$	 &3.3279E-02  &  1.0597 & 4.4027E-03& 2.0147& 6.3922E-03  &2.0133&3.2254E-03&2.1660  \tabularnewline	\midrule 
				$1/32$& 1.6393E-02 &  1.0215 & 1.0989E-03&2.0023& 1.5939E-03  & 2.0038&7.8548E-04&2.0378     \tabularnewline	\midrule 
				$1/64$& 8.1559E-03  & 1.0072  &2.7458E-04 &2.0008&  3.9819E-04 & 2.0010&1.9511E-04& 2.0093   \tabularnewline	\midrule 
				$1/128$	& 4.0720E-03  & 1.0021  &6.8639E-05&2.0001& 9.9529E-05   & 2.0003&4.8698E-05 & 2.0024  \tabularnewline
				
				\bottomrule	
			\end{tabular}
			\footnotesize
		\end{adjustbox}
		\caption{\label{ex1} Relative errors and convergence rates for example 1, structured mesh.}
		
	\end{table}
	
	\begin{table}
		\centering
		\begin{adjustbox}{max width=\textwidth}
			\begin{tabular}{ccccccccc}
				\toprule 
				\multicolumn{9}{c}{Method 1-rotation on the interaction region}\tabularnewline
				\midrule 
				\multirow{2}{*}{$h$}  & \multicolumn{2}{c}{$\|\underline{\sigma}-\underline{\sigma}_h\|_{L^2(\Omega)}$} & \multicolumn{2}{c}{$\|\underline{\sigma}-\mathcal{M}\underline{\sigma}_h\|_{L^2(\Omega)}$}& \multicolumn{2}{c}{$\|\bm{u}-\bm{u}_h\|_{L^2(\Omega)}$}& \multicolumn{2}{c}{$\|\gamma-\gamma_h\|_{L^2(\Omega)}$}\tabularnewline
				\cmidrule{2-9} \cmidrule{3-9} \cmidrule{4-9} \cmidrule{5-9} \cmidrule{6-9} \cmidrule{7-9} \cmidrule{8-9} \cmidrule{9-9} 
				& Error & Rate & Error & Rate &Error & Rate&Error & Rate\tabularnewline
				\midrule 
				$1/4$& 2.2111E-01 &/   & 3.7442E-02& /& 2.5009E-03  & /& 9.4630E-02& /  \tabularnewline\midrule 
				$1/8$&  9.9777E-02& 1.1480  & 7.1351E-03& 2.3917&   9.4662E-04&1.4016&3.3204E-02& 1.5109  \tabularnewline\midrule 
				$1/16$	  &4.8401E-02  & 1.0437   &1.7323E-03 &2.0422&2.9016E-04   &1.7059&1.1610E-02& 1.5160  \tabularnewline\midrule 
				$1/32$& 2.4226E-02 &0.9985  &4.9775E-04 &1.7992& 7.9064E-05   & 1.8758&4.0702E-03  &  1.5122 \tabularnewline\midrule 
				$1/64$ & 1.2184E-02 &  0.9916 & 1.5365E-04 &1.6958&  2.0493E-05  &1.9479&1.4311E-03  & 1.5080 \tabularnewline\midrule 
				$1/128$	&6.1189E-03  &0.9936   &4.9127E-05  &1.6451 &   5.2025E-06& 1.9779& 5.0423E-04 & 1.5050 \tabularnewline	\midrule 
				$1/256$	& 3.0673E-03 &  0.9963 &1.6135E-05 & 1.6063& 1.3012E-06 & 1.9994& 1.7791E-04 &1.5029  \tabularnewline
				\toprule 
				\multicolumn{9}{c}{Method 2- rotation on the macro-element}\tabularnewline
				\midrule 
				\multirow{2}{*}{$h$}  & \multicolumn{2}{c}{$\|\underline{\sigma}-\underline{\sigma}_h\|_{L^2(\Omega)}$} & \multicolumn{2}{c}{$\|\underline{\sigma}-\mathcal{M}\underline{\sigma}_h\|_{L^2(\Omega)}$}& \multicolumn{2}{c}{$\|\bm{u}-\bm{u}_h\|_{L^2(\Omega)}$}& \multicolumn{2}{c}{$\|\gamma-\gamma_h\|_{L^2(\Omega)}$}\tabularnewline
				\cmidrule{2-9} \cmidrule{3-9} \cmidrule{4-9} \cmidrule{5-9} \cmidrule{6-9} \cmidrule{7-9} \cmidrule{8-9} \cmidrule{9-9}
				& Error & Rate & Error & Rate &Error & Rate&Error & Rate\tabularnewline
				\midrule 
				$1/4$ &1.1010E-01 & /  &1.8332E-03 & /&  4.6044E-05  &/& 4.7322E-04 & /\tabularnewline\midrule
				$1/8$  &   5.4722E-02&   1.0086& 4.8802E-04 & 1.9094&  1.1665E-05  &1.9808&1.3065E-04  & 1.8568\tabularnewline\midrule 
				$1/16$	 &  2.7326E-02 &  1.0018 &  1.3006E-04&1.9078 & 3.1004E-06  & 1.9117& 3.4599E-05& 1.9169  \tabularnewline\midrule 
				$1/32$&  1.3659E-02 &  1.0004 &3.4524E-05&1.9135 &7.9405E-07  & 1.9652&  9.0687E-06&  1.9318 \tabularnewline\midrule 
				$1/64$ & 6.8288E-03 & 1.0001  & 9.1124E-06&1.9217&   2.0016E-07&1.9881&  2.3667E-06 &  1.9380 \tabularnewline\midrule  
				$1/128$	 &3.4143E-03  &   1.0000&2.3924E-06  &1.9294 &  5.0178E-08 &1.9960&6.2229E-07 &1.9272  \tabularnewline

				\bottomrule
			\end{tabular}
		}
	\end{adjustbox}
	\caption{\label{ex2} Relative errors and convergence rates for example 2, three dimensional tests, structured mesh.}
	\label{}
\end{table}

\begin{table}

\centering
\begin{adjustbox}{max width=\textwidth}
	
	\begin{tabular}{ccccccccc}
		\toprule 
		\multicolumn{9}{c}{Method 1-scaled rotation on the interaction region }\tabularnewline
		\midrule 
		\multirow{2}{*}{$h$}  & \multicolumn{2}{c}{$\|\underline{\sigma}-\underline{\sigma}_h\|_{L^2(\Omega)}$} & \multicolumn{2}{c}{$\|\underline{\sigma}-\mathcal{M}\underline{\sigma}_h\|_{L^2(\Omega)}$}& \multicolumn{2}{c}{$\|\bm{u}-\bm{u}_h\|_{L^2(\Omega)}$}& \multicolumn{2}{c}{$\| \tilde{\gamma}-\tilde{\gamma}_h\|_{L^2(\Omega)}$}\tabularnewline
		\cmidrule{2-9} \cmidrule{3-9} \cmidrule{4-9} \cmidrule{5-9} \cmidrule{6-9} \cmidrule{7-9} \cmidrule{8-9} \cmidrule{9-9} 
		& Error & Rate & Error & Rate &Error & Rate&Error & Rate\tabularnewline
		\midrule 
		$1/6$&  4.3083E-01 &   /& 2.4907E-01&/&  3.4472E-01 & /&5.9065E-01  &  /\tabularnewline\midrule 
		$1/12$	 &2.0504E-01 & 1.0712 &7.2653E-02 &1.7775 &9.0011E-02   &1.9373& 3.0859E-01& 0.9366  \tabularnewline\midrule 
		$1/24$& 1.0119E-01 &  1.0188 &2.1522E-02 &1.7552&  2.4867E-02 & 1.8559&1.1538E-01& 1.4193   \tabularnewline\midrule 
		$1/48$& 5.0426E-02 &  1.0048 &6.2383E-03  &1.7866& 6.5215E-03  & 1.9310&3.8841E-02 &  1.5707 \tabularnewline\midrule 
		$1/96$	& 2.5183E-02 &1.0017   & 1.8072E-03& 1.7874& 1.6566E-03  & 1.9770&1.3056E-02& 1.5729  \tabularnewline				
		\toprule 
		\multicolumn{9}{c}{Method 2- rotation on the macro-element}\tabularnewline
		\midrule 
		\multirow{2}{*}{$h$}  & \multicolumn{2}{c}{$\|\underline{\sigma}-\underline{\sigma}_h\|_{L^2(\Omega)}$} & \multicolumn{2}{c}{$\|\underline{\sigma}-\mathcal{M}\underline{\sigma}_h\|_{L^2(\Omega)}$}& \multicolumn{2}{c}{$\|\bm{u}-\bm{u}_h\|_{L^2(\Omega)}$}& \multicolumn{2}{c}{$\|\gamma-\gamma_h\|_{L^2(\Omega)}$}\tabularnewline
		\cmidrule{2-9} \cmidrule{3-9} \cmidrule{4-9} \cmidrule{5-9} \cmidrule{6-9} \cmidrule{7-9} \cmidrule{8-9} \cmidrule{9-9}
		& Error & Rate & Error & Rate &Error & Rate&Error & Rate\tabularnewline
		\midrule 
		$1/6$& 3.8757E-01& /& 2.2966E-01 &/ &2.8285E-01   &/&2.8291E-01  & /  \tabularnewline\midrule 
		$1/12$	& 1.6938E-01& 1.1942& 5.7321E-02& 2.0024&6.7366E-02   &2.0699&5.3597E-02&  2.4001   \tabularnewline\midrule 
		$1/24$&8.2636E-02 &1.0354&1.4342E-02 & 1.9988 & 1.6690E-02  &2.0130&1.2312E-02  &    2.1221 \tabularnewline\midrule 
		$1/48$&4.1202E-02  &1.0041&  3.5902E-03  &1.9981& 4.1680E-03  & 2.0016&3.1044E-03    & 1.9877 \tabularnewline\midrule 
		$1/96$	& 2.0604E-02&0.9998&8.9868E-04 &1.9982& 1.0422E-03  &1.9997& 8.1350E-04&      1.9321 \tabularnewline
		
		\bottomrule
	\end{tabular}
\end{adjustbox}
\caption{\label{ex3} Relative errors and convergence rates for example 3, highly heterogeneous material tests, structured mesh.}

\end{table}

\begin{table}
\centering
\begin{adjustbox}{max width=\textwidth}
	
	\begin{tabular}{ccccccccc}
		\toprule 
		\multicolumn{9}{c}{Method 1-rotation on the interaction region}\tabularnewline
		\midrule 
		\multirow{2}{*}{$h$}  & \multicolumn{2}{c}{$\|\underline{\sigma}-\underline{\sigma}_h\|_{L^2(\Omega)}$} & \multicolumn{2}{c}{$\|\underline{\sigma}-\mathcal{M}\underline{\sigma}_h\|_{L^2(\Omega)}$}& \multicolumn{2}{c}{$\|\bm{u}-\bm{u}_h\|_{L^2(\Omega)}$}& \multicolumn{2}{c}{$\|\gamma-\gamma_h\|_{L^2(\Omega)}$}\tabularnewline
		\cmidrule{2-9} \cmidrule{3-9} \cmidrule{4-9} \cmidrule{5-9} \cmidrule{6-9} \cmidrule{7-9} \cmidrule{8-9} \cmidrule{9-9} 
		& Error & Rate & Error & Rate &Error & Rate&Error & Rate\tabularnewline
		\midrule 
		$1/4$& 3.4578E-01  & /  & 6.2566E-02&/ & 7.7210E-02  & /& 1.9954E-01& /\tabularnewline\midrule
		$1/8$& 1.6957E-01 & 1.0280  &2.2109E-02& 1.5007& 2.0867E-02   & 1.8876& 7.8194E-02& 1.3515  \tabularnewline\midrule 
		$1/16$	&8.9611E-02  &0.9201   &7.0519E-03 &1.6485&5.3078E-03   & 1.9750&2.6940E-02&  1.5373 \tabularnewline\midrule 
		$1/32$& 4.6001E-02 &  0.9620 & 1.9090E-03&1.8852& 1.3310E-03  & 1.9956&9.1864E-03&   1.5522 \tabularnewline\midrule 
		$1/64$ & 2.3153E-02 & 0.9905  & 4.8827E-04&1.9671& 3.3295E-04  &1.9991&3.1666E-03& 1.5366   \tabularnewline\midrule 
		$1/128$	& 1.1596E-02 & 0.9976  &1.2287E-04 & 1.9905&   8.3247E-05 & 1.9998&1.1030E-03&1.5215   \tabularnewline				
		\toprule 
		\multicolumn{9}{c}{Method 2- rotation on the macro-element}\tabularnewline
		\midrule 
		\multirow{2}{*}{$h$}  & \multicolumn{2}{c}{$\|\underline{\sigma}-\underline{\sigma}_h\|_{L^2(\Omega)}$} & \multicolumn{2}{c}{$\|\underline{\sigma}-\mathcal{M}\underline{\sigma}_h\|_{L^2(\Omega)}$}& \multicolumn{2}{c}{$\|\bm{u}-\bm{u}_h\|_{L^2(\Omega)}$}& \multicolumn{2}{c}{$\|\gamma-\gamma_h\|_{L^2(\Omega)}$}\tabularnewline
		\cmidrule{2-9} \cmidrule{3-9} \cmidrule{4-9} \cmidrule{5-9} \cmidrule{6-9} \cmidrule{7-9} \cmidrule{8-9} \cmidrule{9-9}
		& Error & Rate & Error & Rate &Error & Rate&Error & Rate\tabularnewline
		\midrule 
		$1/4$& 3.4495E-01 & /  & 3.6450E-02& /& 7.3956E-02  & /&2.4462E-02& /\tabularnewline\midrule 
		$1/8$& 1.7079E-01   &   1.0142& 7.6524E-03& 2.2519& 1.9286E-02  & 1.9391&6.4972E-03&1.9127   \tabularnewline\midrule 
		$1/16$& 9.0944E-02  &  0.9092 & 1.8491E-03&2.0491	& 4.8739E-03  &1.9844 &1.6847E-03&1.9473  \tabularnewline\midrule 
		$1/32$ &4.6223E-02  &0.9764   &4.6841E-04&1.9810& 1.2210E-03   &1.9970& 4.2476E-04&1.9878    \tabularnewline\midrule 
		$1/64$ & 2.3182E-02 & 0.9956  &1.1794E-04 &1.9897& 3.0536E-04  &1.9995&1.0636E-04 &   1.9977 \tabularnewline\midrule 
		$1/128$	 & 1.1600E-02 & 0.9989  &2.9552E-05& 1.9967&7.6345E-05   &1.9999& 2.6600E-05&1.9995  \tabularnewline
		
		\bottomrule
	\end{tabular}
	
\end{adjustbox}
\caption{ \label{ex4}Relative errors and convergence rates for example 4,  looking free tests, $\mu=1$ and $\lambda=10^6$, structured mesh.}
\end{table}

\begin{figure}
\centering
\subfigure[ Parallelogram mesh, $h=1/32$.]{
	\includegraphics[trim={2cm 0.3cm 1cm 0.2cm},clip,width=2in]{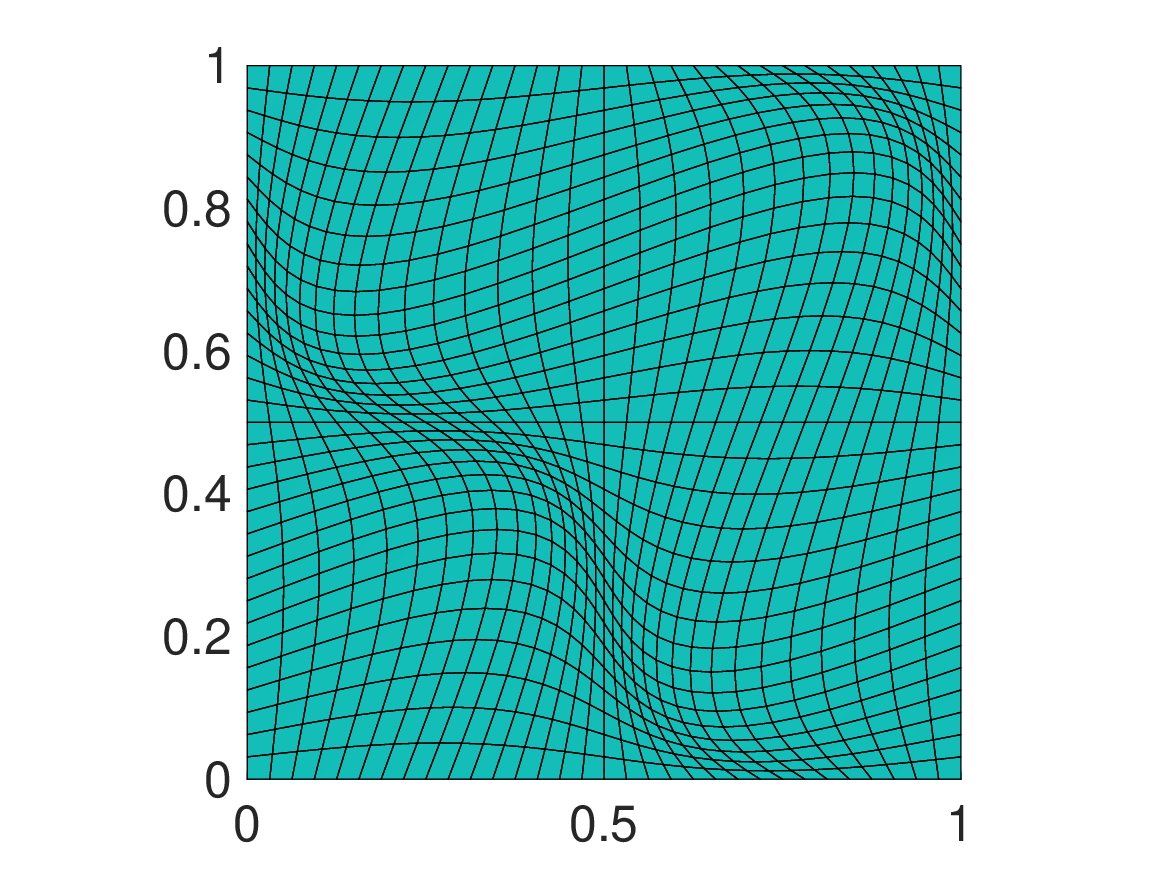}\label{fig:grid_quad}}	
\subfigure[Smooth quadrilateral mesh, $h=1/32$.]{
	\includegraphics[trim={2cm 0.3cm 1cm 0.2cm},clip,width=2in]{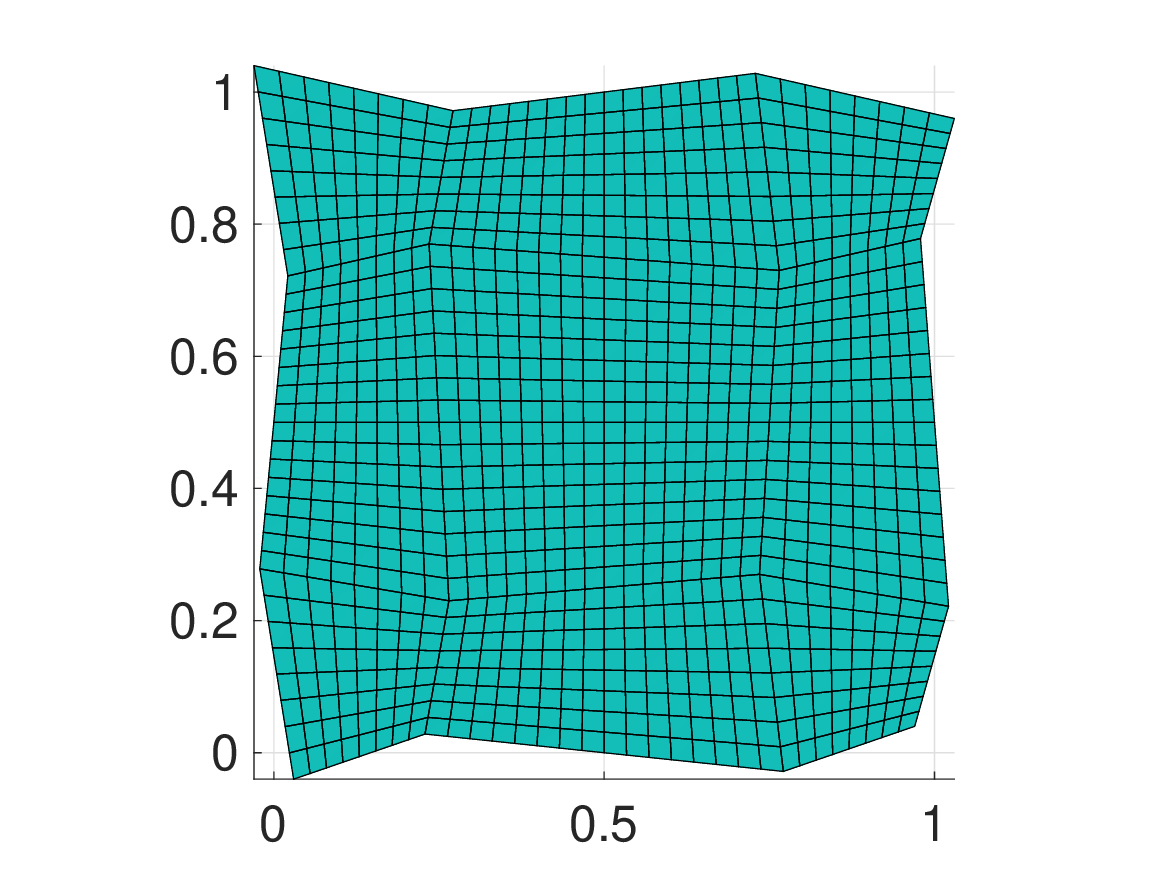}\label{fig:grid_smooth}}	
\subfigure[Randomly perturbed grids  $h^2$, $h=1/32$.]{
	\includegraphics[trim={2cm 0.3cm 1cm 0.2cm},clip,width=2in]{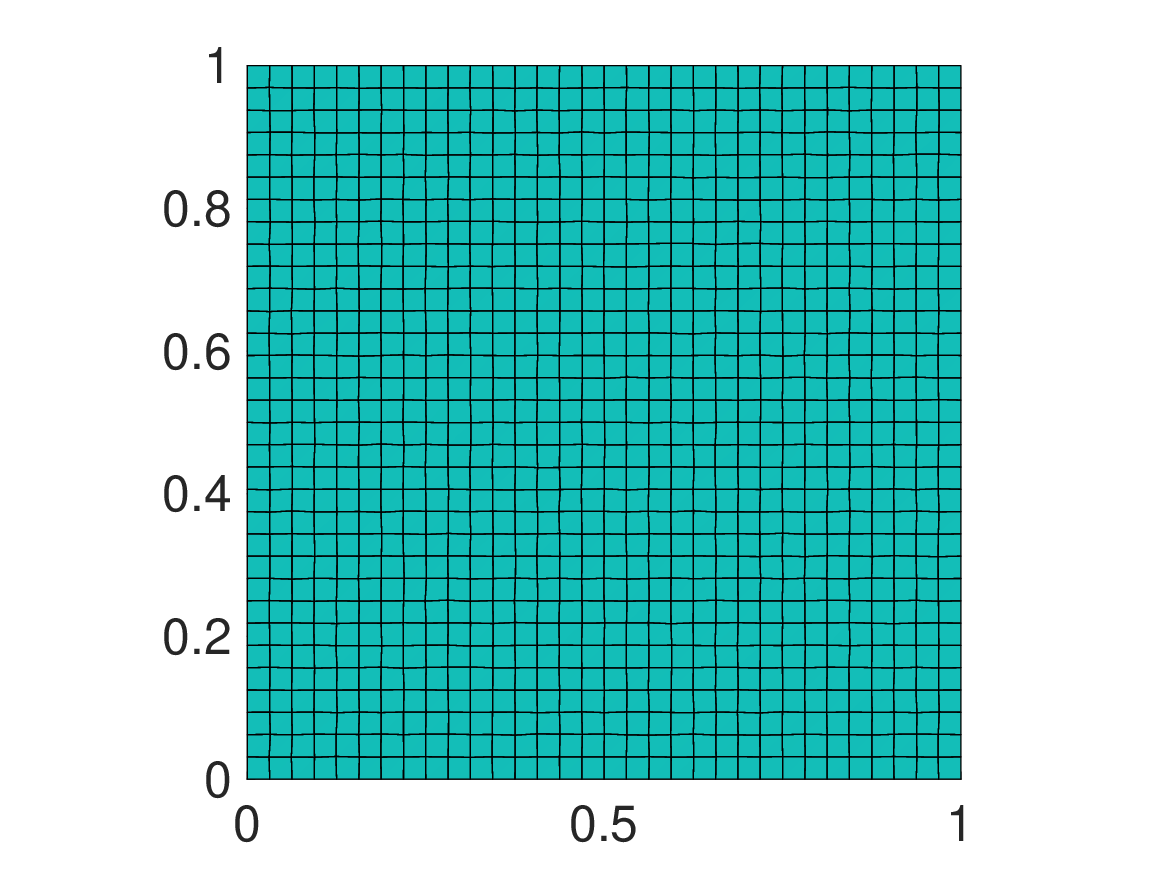}\label{fig:grid_h2}}								
\caption{Examples of unstructured meshes.}
\label{fig:griduns}
\end{figure}

\hong{\subsection{Example 5: 2D unstructured grid cases}}
We next study convergence of the proposed methods on three different types of unstructured grids, we only report results of rotation defined over the interaction region (Method 1). 
For the first test, we first partition the unit square into a $4\times 4$ square mesh with $h=\frac{1}{4}$, and then move the mesh points by applying a map 
\begin{equation*}
x=\hat{x}+0.03\cos(3\pi \hat{x})\cos(3\pi \hat{y}),\quad 	y=\hat{y}-0.04\cos(3\pi \hat{x})\cos(3\pi \hat{y})
\end{equation*}
which yields a deformed computational domain with $4\times 4$ quadrilateral grid. We then apply a uniform refinement of this $4\times 4$ quadrilateral grid to generate a sequence of grids (\hong{see an example in Figure \ref{fig:grid_quad}}), the resulting sequence of mesh satisfies the $h^2$-parallelogram property. Convergence rates on this grid are shown in Table 
\ref{ex1_para}, at least first-order convergence for all variables in their respective
norms are observed, there are no significant degeneration of the convergence rate. For the second unstructured grid test, we compute solution on a sequence of smooth quadrilateral meshes which is obtained by applying a smooth map 
\begin{equation*}
x=\hat{x}+0.1\sin(2\pi \hat{x})\sin(2\pi \hat{y}),\quad 	y=\hat{y}+0.1\sin(2\pi \hat{x})\sin(2\pi \hat{y})
\end{equation*}
to a uniformly refined square mesh (\hong{see an example in Figure \ref{fig:grid_smooth}}). We reported errors and convergence rates on this meshes in Table \ref{ex1_sm},  similar convergence behaviors as parallelogram mesh are found.
Displacement computed on this type of mesh are depicted in right panel of Figure \ref{fig:sol}. 

We also test our method on sequences of quadrilateral grids
generated by random perturbation of uniform grids (\hong{see an example in  Figure \ref{fig:grid_h2}}). Specifically, at each refinement level, the vertices are randomly perturbed within a circle whose radius is $h^2$. 
Results for these meshes are reported in Table \ref{ex1_h2}, again, the proposed method delivers superconvergence of displacement and at least first-order convergence of other variables.
\begin{table}
\centering
\begin{adjustbox}{max width=\textwidth}
	
	\begin{tabular}{ccccccccc}
		\toprule 
		
		\multirow{2}{*}{$h$}  & \multicolumn{2}{c}{$\|\underline{\sigma}-\underline{\sigma}_h\|_{L^2(\Omega)}$} & \multicolumn{2}{c}{$\|\underline{\sigma}-\mathcal{M}\underline{\sigma}_h\|_{L^2(\Omega)}$}& \multicolumn{2}{c}{$\|\bm{u}-\bm{u}_h\|_{L^2(\Omega)}$}& \multicolumn{2}{c}{$\|\gamma-\gamma_h\|_{L^2(\Omega)}$}\tabularnewline
		\cmidrule{2-9} \cmidrule{3-9} \cmidrule{4-9} \cmidrule{5-9} \cmidrule{6-9} \cmidrule{7-9} \cmidrule{8-9} \cmidrule{9-9} 
		& Error & Rate & Error & Rate &Error & Rate&Error & Rate\tabularnewline
		\midrule

		$1/4$& 2.1107E-01&/   &9.2463E-02& /& 1.3043E-01&/&   1.5986E-01	 &/\tabularnewline\midrule
		$1/8$ &1.0759E-01& 0.9722  &3.2100E-02&1.5263&3.1790E-02 &2.0366&6.3990E-02  &1.3209\tabularnewline\midrule
		$1/16$	 &  5.3408E-02& 1.0104 &8.8784E-03& 1.8542& 8.0459E-03 &1.9822&  2.7066E-02  &1.2414 \tabularnewline\midrule
		$1/32$ &2.6938E-02& 0.9874 &2.5908E-03&1.7769& 2.0474E-03&1.9745&9.9050E-03 & 1.4503\tabularnewline\midrule
		$1/64$ &1.3501E-02&  0.9966  &7.8447E-04& 1.7236& 5.1790E-04&1.9830& 3.3807E-03 &1.5508 \tabularnewline\midrule
		$1/128$	 &6.7604E-03& 0.9979 &2.4824E-04& 1.6600&1.3003E-04 &1.9938& 1.1457E-03 &1.5611 \tabularnewline
		\bottomrule
	\end{tabular}
	
\end{adjustbox}
\caption{	\label{ex1_para} Relative errors and convergence rates for Example 1, parallelogram mesh, Method 1- rotation on the interaction region.}

\end{table}


\begin{table}
\centering
\begin{adjustbox}{max width=\textwidth}
	
	\begin{tabular}{ccccccccc}
		\toprule 
		
		\multirow{2}{*}{$h$}  & \multicolumn{2}{c}{$\|\underline{\sigma}-\underline{\sigma}_h\|_{L^2(\Omega)}$} & \multicolumn{2}{c}{$\|\underline{\sigma}-\mathcal{M}\underline{\sigma}_h\|_{L^2(\Omega)}$}& \multicolumn{2}{c}{$\|\bm{u}-\bm{u}_h\|_{L^2(\Omega)}$}& \multicolumn{2}{c}{$\|\gamma-\gamma_h\|_{L^2(\Omega)}$}\tabularnewline
		\cmidrule{2-9} \cmidrule{3-9} \cmidrule{4-9} \cmidrule{5-9} \cmidrule{6-9} \cmidrule{7-9} \cmidrule{8-9} \cmidrule{9-9} 
		& Error & Rate & Error & Rate &Error & Rate&Error & Rate\tabularnewline
		\midrule

		$1/4$&2.2844E-01&/   &1.2264E-01& /& 1.4740E-01&/&   2.4499E-01&/\tabularnewline\midrule
		$1/8$ & 1.3232E-01& 0.7878 &4.8392E-02&1.3416&4.7768E-02&1.6256&   1.3381E-01 &0.8725 \tabularnewline\midrule
		$1/16$	 &6.7906E-02&   0.9624&1.6935E-02& 1.5148& 1.4576E-02&1.7124&  6.6473E-02& 1.0093\tabularnewline\midrule
		$1/32$ & 3.4512E-02&  0.9764&5.3464E-03&1.6634&4.0691E-03 & 1.8408&   2.3594E-02 &1.4943 \tabularnewline\midrule
		$1/64$ & 1.7404E-02&  0.9877& 1.6989E-03&1.6540 & 1.0603E-03&1.9402& 7.5525E-03   &1.6434 \tabularnewline\midrule
		$1/128$	 &8.7378E-03& 0.9941&5.6783E-04&1.5811&2.6858E-04&1.9810&  2.4374E-03 & 1.6316\tabularnewline
		\bottomrule
	\end{tabular}
	
\end{adjustbox}
\caption{	\label{ex1_sm} Relative errors and convergence rates for example 1, smooth quadrilateral mesh, Method 1-rotation on the interaction region.}

\end{table}

\begin{table}
\centering
\begin{adjustbox}{max width=\textwidth}
	
	\begin{tabular}{ccccccccc}
		\toprule 
		
		\multirow{2}{*}{$h$}  & \multicolumn{2}{c}{$\|\underline{\sigma}-\underline{\sigma}_h\|_{L^2(\Omega)}$} & \multicolumn{2}{c}{$\|\underline{\sigma}-\mathcal{M}\underline{\sigma}_h\|_{L^2(\Omega)}$}& \multicolumn{2}{c}{$\|\bm{u}-\bm{u}_h\|_{L^2(\Omega)}$}& \multicolumn{2}{c}{$\|\gamma-\gamma_h\|_{L^2(\Omega)}$}\tabularnewline
		\cmidrule{2-9} \cmidrule{3-9} \cmidrule{4-9} \cmidrule{5-9} \cmidrule{6-9} \cmidrule{7-9} \cmidrule{8-9} \cmidrule{9-9} 
		& Error & Rate & Error & Rate &Error & Rate&Error & Rate\tabularnewline
		\midrule

		$1/4$&1.9552E-01 &/   &9.8976E-02& /& 1.3217E-01&/&  2.1107E-01 &/\tabularnewline\midrule
		$1/8$ &1.0162E-01&   0.9441&2.7386E-02&1.8536 & 2.8569E-02& 2.2099& 4.5648E-02& 2.2091\tabularnewline\midrule
		$1/16$	 & 4.9652E-02 &  1.0333 & 7.5703E-03& 1.8550& 7.1306E-03&2.0024& 1.3117E-02  &1.7991 \tabularnewline\midrule
		$1/32$ &2.4909E-02& 0.9952  &2.1222E-03& 1.8348& 1.7565E-03&2.0213&3.6408E-03  &1.8491 \tabularnewline\midrule
		$1/64$ &1.2447E-02&  1.0009 & 6.7770E-04& 1.6468& 4.3809E-04 &2.0034&  1.2704E-03 &1.5190 \tabularnewline\midrule
		$1/128$	 &6.2262E-03&  0.9994 &2.8947E-04& 1.2272&  1.0955E-04& 1.9996& 5.3802E-04  & 1.2396\tabularnewline
		\bottomrule
	\end{tabular}
	
\end{adjustbox}
\caption{\label{ex1_h2} Relative errors and convergence rates for example 1, randomly perturbed grids  $h^2$, Method 1-rotation on the interaction region.}

\end{table}

\bibliographystyle{plain}
\bibliography{reference}

\begin{thebibliography}{10}

\bibitem{Aavatsmark02}
I.~Aavatsmark.
\newblock An introduction to multipoint flux approximations for quadrilateral
  grids.
\newblock {\em Comput. Geosci.}, 6:405--432, 2002.

\bibitem{Aavatsmark98}
I.~Aavatsmark, T.~Barkve, O.~B\o{}e, and T.~Mannseth.
\newblock Discretization on unstructured grids for inhomogeneous, anisotropic
  media. part {II}: {D}iscussion and numerical results.
\newblock {\em SIAM J. Sci. Comput.}, 19(5):1717--1736, 1998.

\bibitem{Agelas08}
L.~Agelas and R.~Masson.
\newblock Convergence of the finite volume {MPFA} {O} scheme for heterogeneous
  anisotropic diffusion problems on general meshes.
\newblock {\em C. R. Math.}, 346(17):1007--1012, 2008.

\bibitem{Ambartsumyan20simplicial}
I.~Ambartsumyan, E.~Khattatov, J.~M. Nordbotten, and I.~Yotov.
\newblock A multipoint stress mixed finite element method for elasticity on
  simplicial grids.
\newblock {\em SIAM J. Numer. Anal.}, 58(1):630--656, 2020.

\bibitem{Ambartsumyan21}
I.~Ambartsumyan, E.~Khattatov, J.~M. Nordbotten, and I.~Yotov.
\newblock A multipoint stress mixed finite element method for elasticity on
  quadrilateral grids.
\newblock {\em Numer. Meth. Part Differ. Equ.}, 37(3):1886--1915, 2021.

\bibitem{Arnold14}
D.~Arnold, G.~Awanou, and R.~Winther.
\newblock Nonconforming tetrahedral mixed finite elements for elasticity.
\newblock {\em Math. ModelsMethods Appl. Sci.}, 24(04):783--796, 2014.

\bibitem{arnold2007mixed}
D.~Arnold, R.~Falk, and R.~Winther.
\newblock Mixed finite element methods for linear elasticity with weakly
  imposed symmetry.
\newblock {\em Math. Comp.}, 76(260):1699--1723, 2007.

\bibitem{Arnold84peer}
D.~N. Arnold, F.~Brezzi, and J.~Douglas.
\newblock Peers: a new mixed finite element for plane elasticity.
\newblock {\em Japan J. Appl. Math.}, 1:347--367, 1984.

\bibitem{Arnold02}
D.~N. Arnold and R.~Winther.
\newblock Mixed finite elements for elasticity.
\newblock {\em Numer. Math.}, 92:401--419, 2002.

\bibitem{Arnold84}
D.N. Arnold, J.~Douglas, and C.~P. Gupta.
\newblock A family of higher order mixed finite element methods for plane
  elasticity.
\newblock {\em Numer. Math.}, 45:1--22, 1984.

\bibitem{Boffi09}
D.~Boffi, F.~Brezzi, and M.~Fortin.
\newblock Reduced symmetry elements in linear elasticity.
\newblock {\em Commun. Pure Appl. Anal}, 8(1):95--121, 2009.

\bibitem{Brenner03}
S.~C. Brenner.
\newblock Poincaré--{F}riedrichs inequalities for piecewise ${H}^1$ functions.
\newblock {\em SIAM J. Numer. Anal.}, 41(1):306--324, 2003.

\bibitem{Brezzi91}
F.~Brezzi and M.~Fortin.
\newblock {\em Mixed and Hybrid Finite Element Methods}.
\newblock Springer-Verlag, Berlin, 1991.

\bibitem{Edwards02}
M.~G. Edwards.
\newblock Unstructured, control-volume distributed, full-tensor finite-volume
  schemes with flow based grids.
\newblock {\em Comput. Geosci.}, 6:433--452, 2002.

\bibitem{Edwards98}
M.~G. Edwards and C.~F. Rogers.
\newblock Finite volume discretization with imposed flux continuity for the
  general tensor pressure equation.
\newblock {\em Comput. Geosc.}, 2:259--290, 1998.

\bibitem{Ewing99}
R.~E. Ewing, M.~M. Liu, and J.~Wang.
\newblock Superconvergence of mixed finite element approximations over
  quadrilaterals.
\newblock {\em SIAM J. Numer. Anal.}, 36(3):772--787, 1999.

\bibitem{Girault79}
V.~Girault and P.-A. Raviart.
\newblock {\em Finite Element Approximation of the Navier-Stokes Equations}.
\newblock Springer Berlin, Heidelberg, 1979.

\bibitem{Gopalakrishnan11}
J.~Gopalakrishnan and J.~Guzm{\'a}n.
\newblock Symmetric nonconforming mixed finite elements for linear elasticity.
\newblock {\em SIAM J. Numer. Anal.}, 49(4):1504--1520, 2011.

\bibitem{Grisvard11}
P.~Grisvard.
\newblock {\em Elliptic Problems in Nonsmooth Domains}.
\newblock Society for Industrial and Applied Mathematics, 2011.

\bibitem{HornJohnson2012}
R.~A. Horn and C.~R. Johnson.
\newblock {\em Matrix Analysis}.
\newblock Cambridge University Press, 2 edition, 2012.

\bibitem{Ingram10}
R.~Ingram, M.~F. Wheeler, and I.~Yotov.
\newblock A multipoint flux mixed finite element method on hexahedra.
\newblock {\em SIAM J. Numer. Anal.}, 48(4):1281--1312, 2010.

\bibitem{Kechkar92}
N.~Kechkar and D.~Silvester.
\newblock Analysis of locally stabilized mixed finite element methods for the
  {S}tokes problem.
\newblock {\em Math. Comp.}, 58(197):1--10, 1992.

\bibitem{Klausen06}
R.~A. Klausen and R.~Winther.
\newblock Convergence of multipoint flux approximations on quadrilateral grids.
\newblock {\em Numer. Meth. Part Differ. Equ.}, 22(6):1438--1454, 2006.

\bibitem{klausen2006robust}
R.~A. Klausen and R.~Winther.
\newblock Robust convergence of multi point flux approximation on rough grids.
\newblock {\em Numer. Math.}, 104:317--337, 2006.

\bibitem{PechsteinTDNNS18}
A.~S. Pechstein and J.~Sch{\"o}berl.
\newblock An analysis of the {TDNNS} method using natural norms.
\newblock {\em Numer. Math.}, 139:93--120, 2018.

\bibitem{Qiu18}
W.~Qiu, J.~Shen, and K.~Shi.
\newblock An {HDG} method for linear elasticity with strong symmetric stresses.
\newblock {\em Math. Comp.}, 87:69--93, 2018.

\bibitem{wheeler2012multipoint}
M.~F. Wheeler, G.~Xue, and I.~Yotov.
\newblock A multipoint flux mixed finite element method on distorted
  quadrilaterals and hexahedra.
\newblock {\em Numer. Math.}, 121:165--204, 2012.

\bibitem{Wheeler06}
M.~F. Wheeler and I.~Yotov.
\newblock A multipoint flux mixed finite element method.
\newblock {\em SIAM J. Numer. Anal.}, 44:2082--2106, 2006.

\bibitem{Zhao20}
L.~Zhao and E.-J. Park.
\newblock A staggered cell-centered {DG} method for linear elasticity on
  polygonal meshes.
\newblock {\em SIAM J. Sci. Comput.}, 42(4):A2158--A2181, 2020.

\bibitem{Zhao19}
L.~Zhao, E.-J. Park, and D.~w.~Shin.
\newblock A staggered {DG} method of minimal dimension for the {S}tokes
  equations on general meshes.
\newblock {\em Comput. Methods Appl. Mech. Engrg.}, 345:854--875, 2019.

\end{thebibliography}

\end{document}